\documentclass[aop]{imsart}
\usepackage[utf8]{inputenc}
\usepackage{enumerate}
\usepackage{lmodern}
\usepackage[T1]{fontenc}
\usepackage{verbatim}

\usepackage{textcomp}

\usepackage[english]{babel}
\usepackage[pdftex]{hyperref}

\usepackage[pdftex]{color,graphicx}

\usepackage{amsmath,amsfonts,amssymb,amsthm,mathrsfs}

\newcommand{\e}{\mbox{e}\hspace{-.34em}\mbox{\i}}

\newtheorem{theorem}{Theorem}
\newtheorem{definition}{Definition}[section]
\newtheorem{lemma}[theorem]{Lemma}

\newtheorem{corollary}[theorem]{Corollary}
\newtheorem{proposition}[theorem]{Proposition}

\def\pp{{\mathcal P}}

\def\t{{\mathcal T}}
\def\cc{\mathcal{C}}
\def\ve{\varepsilon}
\def\wt{\widetilde}

\def\f{\mathcal{F}}
\def\g{\mathcal{G}}
\def\h{\mathcal{H}}

\def\QQ{\mathcal{Q}}

\def\e{\mathcal{E}}
\def\la{\longrightarrow}

\def\E{{\mathbb E}}
\def\P{{\mathbb P}}

\def\Z{{\mathbb Z}}

\def\Q{{\mathbb Q}}
\def\F{{\mathbb F}}

\def\C{{\mathbb C}}

\def\dg{\mathrm{d}_{\mathrm{gr}}}

\def\un{\underline}
\def\tr{\mathrm{tr}}
\def\da{\downarrow}

\def\build#1_#2^#3{\mathrel{
\mathop{\kern 0pt#1}\limits_{#2}^{#3}}}
\def\rem{\noindent{\bf Remark. }}

\begin{document}
\begin{frontmatter}

\title{Separating cycles and isoperimetric\\ inequalities
in the uniform infinite \\planar quadrangulation}

\runtitle{Separating cycles and isoperimetric inequalities}

\thankstext{T1}{Supported by the ERC Advanced Grant 740943 {\sc GeoBrown}}

\begin{aug}
\author{\fnms{Jean-Fran\c cois} \snm{Le Gall}\corref{}\ead[label=e1]{jean-francois.legall@math.u-psud.fr}}
\and
\author{\fnms{Thomas} \snm{Leh\'ericy}\ead[label=e2]{thomas.lehericy@math.u-psud.fr}}

\runauthor{J.-F. Le Gall and T. Leh\'ericy}

\affiliation{Universit\'e Paris-Sud}

\address{Institut de Math\'ematique d'Orsay, Universit\'e Paris-Sud, 91405 Orsay, France \\ 
\printead{e1,e2}}

\end{aug}

\begin{abstract}
\ We study geometric properties of the infinite random lattice called the
uniform infinite planar quadrangulation or UIPQ. We establish a precise form of a conjecture of 
Krikun stating that the minimal size of a cycle that separates the
ball of radius $R$ centered at the root vertex from infinity grows linearly in $R$. As a consequence,
we derive certain isoperimetric bounds showing that the 
boundary size of any simply connected set $A$ consisting of a finite union of faces 
of the UIPQ and containing the root vertex is bounded below by a (random) constant times
$|A|^{1/4}(\log|A|)^{-(3/4)-\delta}$, where the volume $|A|$ is the number of faces in $A$.
\end{abstract}



\end{frontmatter}

\section{Introduction}
In the recent years, much work has been devoted to discrete and continuous models 
of random geometry in two dimensions. Two of the most popular discrete models are 
the uniform infinite planar triangulation (or UIPT), which was introduced by Angel and Schramm \cite{Ang,AS}
and in fact motivated much of the subsequent work, and the uniform infinite planar quadrangulation (or UIPQ). In the present work, we concentrate on the UIPQ, although we believe that our
methods can be adapted to give similar results for the UIPT. Roughly speaking, the UIPQ is a
random infinite graph embedded in the plane, such that all faces (connected components of
the complement of edges) are quadrangles, possibly with two edges glued together. See 
Fig.\ref{uipq} below for an illustration of what the UIPQ may look like near its root vertex.
We study certain
geometric properties of the UIPQ, concerning the existence of ``small'' cycles that separate
a large ball centered at the root vertex from infinity, with applications to isoperimetric inequalities.

\begin{figure}[!h]
 \begin{center}
 \includegraphics[width=8cm]{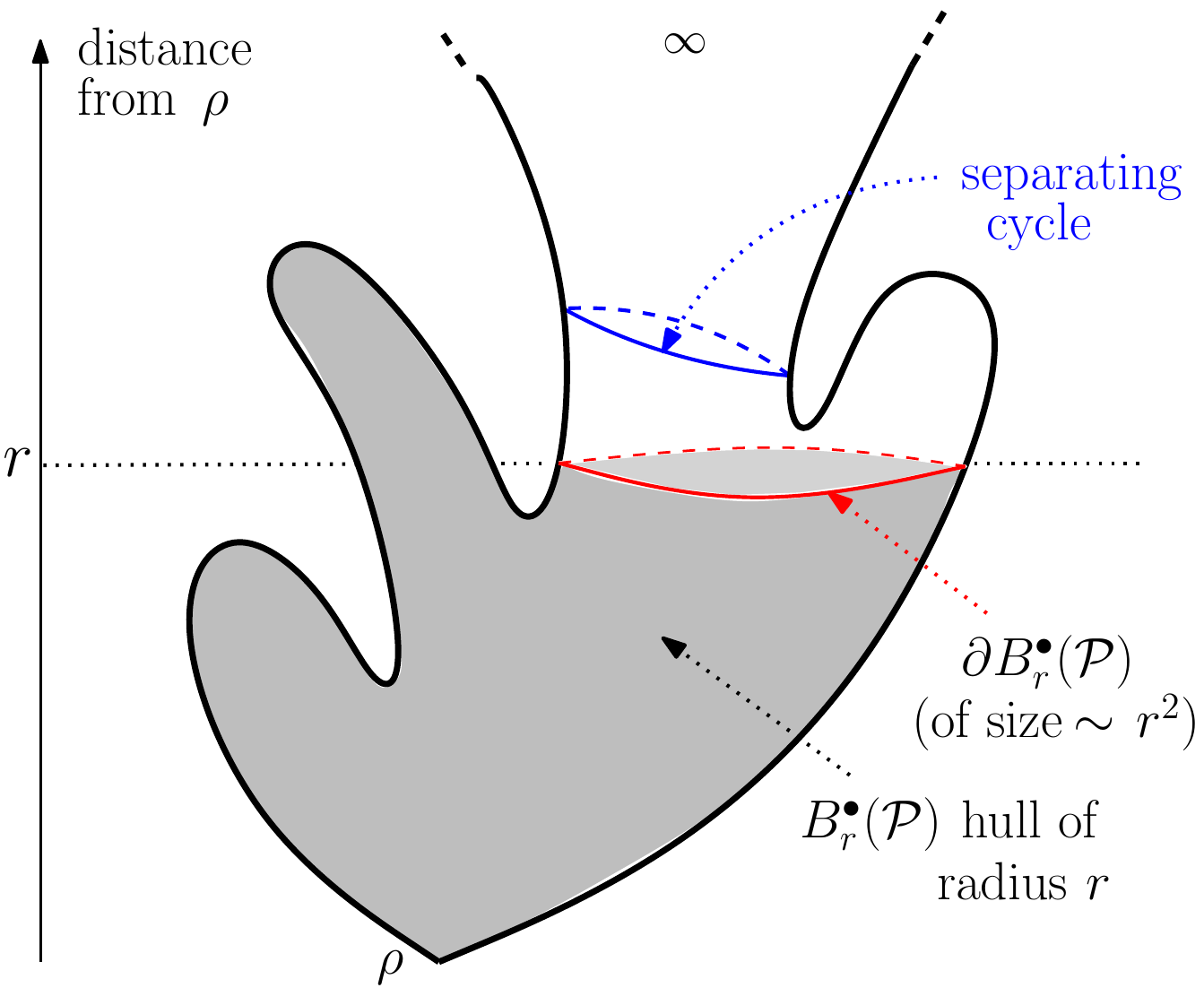}
 \caption{ \label{sep-cyc} A schematic ``cactus'' representation of the UIPQ. The root vertex is denoted by $\rho$
 and the vertical coordinate corresponds to the graph distance from $\rho$. The shaded part
 is the hull $B^\bullet_r(\mathcal{P})$.}
 \vspace{-5mm}
 \end{center}
 \end{figure}

The starting point of our work is a conjecture of Krikun in the paper  \cite{Kr}
which provided the first construction of the UIPQ as the local limit of uniform
planar quadrangulations with a fixed number of faces  (another construction was suggested by
Chassaing and Durhuus \cite{CD}, and the equivalence between the two approaches was
established by M\'enard \cite{Men} --- see also \cite{CMM}
for a third construction). Denote the UIPQ by $\pp$, and, for every integer $r\geq 1$, 
let $B_r(\pp)$ stand for the ball of radius $r$ centered at the root vertex, which is defined 
as the union of all faces that are incident to at least one vertex whose graph distance from
the root is at most $r-1$. The complement of the ball $B_r(\pp)$ is in general 
not connected, but there is a unique unbounded component, whose boundary 
is called the exterior boundary of the ball. The set inside the exterior boundary, 
which may be obtained by filling in the ``bounded holes'' of the ball, is called
the (standard) hull of radius $r$ and will be denoted by $B^\bullet_r(\pp)$. It is known that the size of the exterior boundary, that is, the number of edges
in this boundary, grows like $r^2$ when $r\to\infty$: See \cite{CLG1} for more precise 
asymptotics obtained both for the UIPT and the UIPQ.
On the other hand, Krikun constructed a cycle 
that separates the ball $B_r(\pp)$ from infinity and whose size grows linearly in $r$ when
$r$ is large. Here we say that a cycle $\cc$ made of edges of the UIPQ separates a finite set $A$ of vertices from infinity if 
$A$ does not intersect $\cc$ but
any path from a vertex of $A$ to infinity intersects $\cc$ (see Fig.~\ref{sep-cyc} for
a schematic illustration).  Krikun 
conjectured that the cycle he constructed is essentially the shortest possible, meaning that the minimal size of 
a cycle that separates the ball $B_r(\pp)$ from infinity must be linear in $r$. A weak form
of this conjecture was derived in \cite{plane}, but the results of this paper did not exclude
the possibility that a ball could be separated from infinity by a small cycle lying ``very far
away'' from the ball.

The following theorem provides quantitative estimates that confirm Krikun's conjecture.

\begin{theorem}
\label{kri-conj}
For every integer $R\geq 1$, let $L(R)$ be the smallest length of a cycle separating $B_R(\pp)$ from infinity. 
\begin{itemize}
\item[{\rm(i)}] For every $\delta<2$, there exists a
constant $C_\delta$ such that, for every $R\geq 1$, for every $\ve\in(0,1)$,
$$\P(L(R)\leq \ve R) \leq C_\delta\,\ve^\delta.$$
\item[{\rm(ii)}] There exist constants $C$ and $\lambda>0$ such that, for every $a>0$ and $R\geq 1$,
$$\P(L(R)\geq aR) \leq C\,e^{-\lambda a}.$$
\end{itemize}
\end{theorem}

Part (ii) of the theorem is proved by using the separating cycle introduced by Krikun and 
sharpening the estimates in \cite{Kr}. So the most interesting part of the theorem 
is part (i). We believe that our condition $\delta<2$ is close to optimal, in the sense that, for
$R$ large, $\P(L(R)\leq \ve R)$ should behave like $\ve^2$, possibly up to logarithmic corrections.
At the end of Section \ref{sec:upper}, we provide a short argument showing that 
the probability $\P(L(R)\leq \ve R)$ is bounded below by $\hbox{Const.}\,\ve^3$ when $R$ is large. 

The proof of part (i) relies on a technical estimate which is of
independent interest and that we now present. We first label vertices of the UIPQ 
by their distances from the root vertex, and for every integer $r\geq 1$, we say that a face of
the UIPQ is $r$-simple if the labels of the vertices incident to this face take the three 
values $r-1,r,r+1$ (note that there are faces such that labels of incident vertices take only two values, these faces are called confluent in \cite{CS}). In each $r$-simple face, we draw a ``diagonal'' connecting the
two corners labeled $r$ (these two corners may correspond to the same vertex), and such 
diagonals, which are {\bf not} edges of the UIPQ, are called $r$-diagonals. Then, there
is a ``maximal'' cycle made of $r$-diagonals, which is simple and such that the labels of vertices lying 
in the unbounded component of the complement of this cycle are at least $r+1$. We denote this
maximal cycle by $\cc_r$, and, for $1\leq r< r'$, we define the annulus 
$\cc(r,r')$ as the part of the UIPQ between the cycles $\cc_r$ and $\cc_{r'}$. See Section
\ref{preli} below for more precise definitions. Note that the cycles $\cc_r$ are {\bf not} made of edges
of the UIPQ in contrast with the separating cycles that we consider in Theorem \ref{kri-conj} and 
in the next proposition.

\begin{proposition}
\label{main}
Let $\beta\in(0,3)$. There exists a constant $C'_\beta$ such that, for every integer 
$r\geq 1$ and for every integer $n\geq 1$, the probability that there exists a cycle of the UIPQ
of length smaller than $r$, which is contained in $\cc(nr,(n+2)r)$, does not intersect $\cc_{(n+2)r}$, and disconnects the root vertex
from infinity, is bounded above by $C'_\beta \,n^{-\beta}$.
\end{proposition}


The condition that the cycle does not intersect $\cc_{(n+2)r}$ is included for technical convenience, and could
be removed from the statement. 

The proof of Proposition \ref{main} relies on a ``skeleton decomposition'' of the UIPQ, which is 
already presented in the work of Krikun \cite{Kr}. Our presentation is however different from 
the one in \cite{Kr} and better suited to our purposes. We introduce and use the notion of a truncated quadrangulation, which is
basically a planar map with a boundary, where all faces (distinct from the distinguished one) 
are quadrangles, except for those incident to the boundary, which are triangles (see Section \ref{sec:trunc} for precise definitions). 
The annulus $\cc(r,r')$ can be viewed as a truncated quadrangulation
of the cylinder of height $r'-r$.  Our motivation for introducing truncated quadrangulations comes
from the fact that they allow certain explicit calculations in the UIPQ. For
every integer $r\geq 1$, we define
the ``truncated hull'' of radius $r$ of the UIPQ, which is basically the part of the UIPQ  inside the maximal cycle $\cc_r$
(see Section \ref{sec:trunc-UIPQ} for a precise definition).
This truncated hull is different from the standard hull $B^\bullet_r(\pp)$ introduced above, which had been considered in \cite{CLG0,CLG1} in 
particular, but it is essentially the same object as the hull defined in \cite{Kr}. It turns out that it is possible to
compute the law of the truncated hull in a rather explicit manner (Corollary \ref{law-hull}) and in
particular the law of the perimeter of the hull has a very simple form (Proposition \ref{hull-peri2}). 
These calculations make heavy use of the skeleton decomposition of the UIPQ, and more
generally of the similar decomposition for truncated quadrangulations of the cylinder.
This decomposition involves a forest structure, which was already described by
Krikun \cite[Section 3.2]{Kr} and is similar to the one for triangulations that was discovered in 
\cite{Kr0} and heavily used in the recent work \cite{CLG} dealing with first-passage
percolation on the UIPT. 

Given the forest structure associated with a truncated quadrangulation of the cylinder, the idea
of the proof of Proposition \ref{main} is as follows. One first observes that, with high probability, there 
exist, for some $\delta>0$, more than $n^\delta$ trees with maximal height in the forest coding the annulus $\cc(nr,(n+2)r)$.
For each of these trees, one can find a vertex on the cycle $\cc_{nr}$ (the interior boundary of the annulus)
which is connected to the exterior boundary $\cc_{(n+2)r}$ by a path of length $2r$.
Assuming that there is a cycle of length $r$ in the annulus that disconnects the root vertex 
from infinity, it follows that any two of these particular vertices of $\cc_{nr}$ can be connected by
a path staying in the annulus with length at most $5r$. Results of Curien and Miermont \cite{CM} about
the graph distances between boundary points in infinite quadrangulations with a boundary, show 
that this cannot occur except on a set of small probability.

\medskip
Our lower bounds on the minimal size of separating cycles lead to interesting 
isoperimetric inequalities showing informally that the size of the boundary of a simply connected 
set which is a finite union of faces and contains the root vertex must be at least of the order of
the volume raised to the power $1/4$. The fact that we cannot do better than the power $1/4$
follows from part (ii) in Theorem \ref{kri-conj}, since it is well known \cite{CD,CLG0,CLG1} that the volume of the ball, 
or of the standard hull,
of radius $r$
is of order $r^4$. We refer to \cite[Chapter 6]{LP} for a thorough discussion of
isoperimetric inequalities on infinite graphs. 

Let $\mathcal{K}$ denote the collection of all simply connected compact subsets of the plane that
are  finite unions of faces of the UIPQ (including their boundaries) and contain the
root vertex. For $A\in\mathcal{K}$, the volume of $A$, denoted by $|A|$, is the number 
of faces of the UIPQ contained in $A$, and the boundary size of $A$, denoted by $|\partial A|$, is 
the number of edges in the boundary of $A$.

\begin{theorem}
\label{iso1}
Let $\delta >0$. Then,
$$\inf_{A\in\mathcal{K}} \frac{|\partial A|}{|A|^{\frac{1}{4}}\,(\log|A|)^{-\frac{3}{4}-\delta}}>0\;,\ \hbox{a.s.}$$
\end{theorem}

The exponent $\frac{3}{4}$ in the statement of the theorem is presumably not the optimal one. 
Our method involves estimates for the tail of the distribution of the volume 
of the hull $B^\bullet_r(\pp)$, which are derived from a first moment bound
(Proposition \ref{tail-hull}). We expect that these estimates can be improved,
leading to a better value of the exponent of $\log|A|$ (the results of Riera \cite{Rie} for the Brownian plane suggest that 
one should be able to replace $\frac{3}{4}$ by $\frac{1}{2}$ in the statement of the theorem). On the other hand, one cannot 
hope to replace $|A|^{\frac{1}{4}}\,(\log|A|)^{-\frac{3}{4}-\delta}$ by $|A|^{\frac{1}{4}}$ in the theorem: Simple zero-one
arguments using the separating cycles introduced by Krikun \cite{Kr} (see Section \ref{sec:geodesics}
below) show that there exist sets $A$ such that the ratio $|\partial A|/ |A|^{\frac{1}{4}}$ is arbitrarily
small. 

Still we can state the following proposition. 

\begin{proposition}
\label{iso}
Let $\ve >0$. There exists a constant $c_\ve>0$ such that, for every integer 
$n\geq 1$, the property
$$|\partial A| \geq c_\ve\,n^{1/4}\;,\quad\hbox{for every }A\in\mathcal{K}\hbox{ such that }|A|\geq n\;,$$
holds with probability at least $1-\ve$.
\end{proposition}

As an immediate consequence of Proposition \ref{iso}, we also get that, for every $\ve>0$ and every $M>1$, we can find
a constant $c_{\ve,M}>0$ such that, for every integer $n\geq 1$,
$$\P\Bigg( \inf_{A\in\mathcal{K},\, n\leq |A|\leq Mn} \frac{|\partial A|}{|A|^{\frac{1}{4}}}\geq c_{\ve,M}\Bigg) \geq 1-\ve.$$
Indeed, we just have to take $c_{\ve,M}= c_\ve/M^{1/4}$, with the notation of Proposition \ref{iso}. But, as explained
after the statement of Theorem \ref{iso1}, we cannot lift the constraint $n\leq |A|\leq Mn$ in the 
last display.

The proofs of both Theorem \ref{iso1} and Proposition \ref{iso} rely on Theorem \ref{kri-conj} and on the fact that
the volume of the hull of radius $r$ is of order $r^4$. Assuming that 
$|\partial A|$ is small, then either the root vertex is sufficiently far from $\partial A$, which implies that
a large ball centered at the root vertex is disconnected from infinity by the small cycle $\partial A$
(so that we can use the estimate of Theorem \ref{kri-conj}) or the root vertex is close 
to $\partial A$, but then it follows that the whole set $A$ is contained in the standard hull of 
radius (approximately) equal to the distance from the root vertex to $\partial A$, which implies that
the volume of $A$ cannot be too big (at this point of the argument, in the proof of Theorem \ref{iso1}, we need estimates
for the tail of the distribution of the volume of hulls).

The paper is organized as follows. Section \ref{preli} presents a number of preliminaries,
concerning truncated quadrangulations, their relations with the UIPQ and their skeleton
decompositions, and a number of
related calculations. As mentioned earlier, this section owes a lot to the work of
Krikun \cite{Kr}, and in particular we make use of enumeration results derived in \cite{Kr}.
One additional motivation for deriving the results of Section \ref{preli} in a somewhat 
more precise form than in \cite{Kr} is the fact that we plan to use these results in
a forthcoming work \cite{Leh} on local modifications of distances in the UIPQ, in the spirit of \cite{CLG}.
Proposition \ref{main} is proved in Section \ref{sec:main}, and part (i) of Theorem \ref{kri-conj} easily follows
from this proposition. Section \ref{sec:upper} is devoted to the proof
of part (ii) of Theorem \ref{kri-conj}. This proof relies on the explicit calculation of the distribution 
of the number of trees with maximal height in the forest coding the annulus $\cc(r,r')$ (Proposition \ref{tree-max}). This calculation
is also used to give an easy lower bound for the probability $\P(L(R)\leq \ve R)$.
Section \ref{sec:iso}
contains the proof of Proposition \ref{iso} and Theorem \ref{iso1}. An important ingredient of the proof of Theorem \ref{iso1}
is Proposition \ref{tail-hull}, which provides a first moment bound for the volume of hulls. 
Finally, the Appendix gives
the proof of a technical lemma stated at the end of Section \ref{preli}, which plays
an important role in Section \ref{sec:main}.

\section{Preliminaries}
\label{preli}

\subsection{Truncated quadrangulations}
\label{sec:trunc}

We will consider truncated quadrangulations. Informally, these are quadrangulations
with a simple boundary, where the quadrangles incident to the boundary are replaced
by triangles. A more precise definition is as follows.

\begin{definition}
\label{trunc-quad}
Let $p\geq 1$ be an integer.
A truncated quadrangulation with boundary size $p$ is a planar map $\mathcal{M}$
having a distinguished face $\mathfrak{f}$ with a simple boundary of size $p$ such that:
\begin{enumerate}
\item[$\bullet$] Each edge of the boundary of $\mathfrak{f}$ is incident both to
$\mathfrak{f}$ and to a triangular face of $\mathcal{M}$ and these triangular faces are distinct.
\item[$\bullet$] All faces other than $\mathfrak{f}$ and the triangular faces incident
to the boundary of $\mathfrak{f}$ have degree $4$.
\end{enumerate}
\end{definition}

\begin{figure}[!h]
 \begin{center}
 \includegraphics[width=5cm]{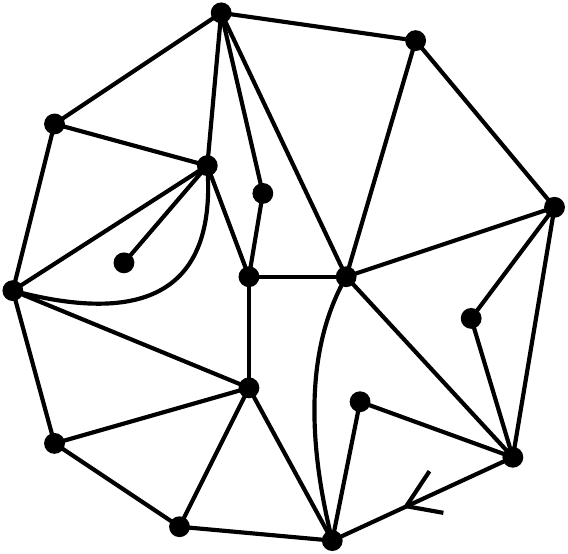}
 \caption{ \label{trunc} A truncated quadrangulation with boundary size $9$, $8$ inner vertices and $16$ inner faces.}
 \vspace{-5mm}
 \end{center}
 \end{figure}

It will be convenient to view truncated quadrangulations as drawn in the plane in such a 
way that the distinguished face is the unbounded face. With this convention, we will 
always assume that a truncated quadrangulation is rooted and, unless otherwise specified, that the root edge 
lies on the boundary of the distinguished face and is oriented clockwise. See Fig.\ref{trunc} for an example.
Faces distinct from the distinguished face are called inner faces, and 
vertices that do not lie on the boundary
of the distinguished face are called inner vertices.

Notice that, when $p\geq 2$, any of the triangular faces  incident to the boundary of $\mathfrak{f}$
must be nondegenerate (i.e. its boundary cannot contain a loop). Furthermore, a simple argument
shows that each of these triangular faces is incident to an inner vertex. The last property clearly also holds if $p=1$. 
Hence a truncated 
quadrangulation with boundary size $p\geq1$ must have at least one inner vertex.

We notice that our truncated quadrangulations with boundary size $p$ are in one-to-one 
correspondence with the ``quadrangulations with a simple boundary'' of size $2p$ 
considered by Krikun \cite{Kr} (starting from the latter, we just ``cut'' the boundary quadrangles
along the appropriate diagonals to get a truncated quadrangulation). If we add an extra vertex $v_*$ inside the face $\mathfrak{f}$, then
draw an edge from each vertex of the boundary of $\mathfrak{f}$ to $v_*$, and finally remove all edges of the boundary of $\mathfrak{f}$, we
get a plane quadrangulation and hence a bipartite graph: In particular, it follows that, if $v$ and $v'$ are two adjacent inner vertices of $\mathcal{M}$, their distances
from the boundary differ by $1$. This observation will be useful later.

For integers $n\geq 1$ and $p\geq 1$, we let 
$\Q^\tr_{n,p}$ be the set of all (rooted) truncated quadrangulations with boundary size $p$
and $n$ inner faces.

We need another definition.

\begin{definition}
\label{trunc-cyl}
Let $h,p,q\geq 1$ be positive integers.
A truncated quadrangulation
of the cylinder of height $h$ with boundary sizes $(p,q)$ is a planar map $\mathcal{Q}$
having two distinguished faces $\mathfrak{f}_b$ and $\mathfrak{f}_t$ such that:
\begin{enumerate}
\item[$\bullet$] The face $\mathfrak{f}_b$ (called the bottom face) has a simple boundary 
of size $p$, which is called the bottom cycle, and the face $\mathfrak{f}_t$ (called the top face) has a simple boundary 
of size $q$, which is called the top cycle. 
\item[$\bullet$] Each edge of the bottom cycle (resp. of the top cycle)  is incident both to
$\mathfrak{f}_b$ (resp. to $\mathfrak{f}_t$) and to a triangular face of $\mathcal{Q}$ and these triangular faces are distinct.
\item[$\bullet$] All faces other than $\mathfrak{f}_b$, $\mathfrak{f}_t$, and the triangular faces incident
to the bottom and top cycles, have degree $4$.
\item[$\bullet$] Every vertex of the top cycle is at graph distance exactly $h$
from the bottom cycle, and every  edge of the top cycle is
incident to a triangular face containing a vertex at graph distance $h-1$ from the bottom cycle. 
\end{enumerate}
\end{definition}

\begin{figure}[!h]
 \begin{center}
 \includegraphics[width=8cm]{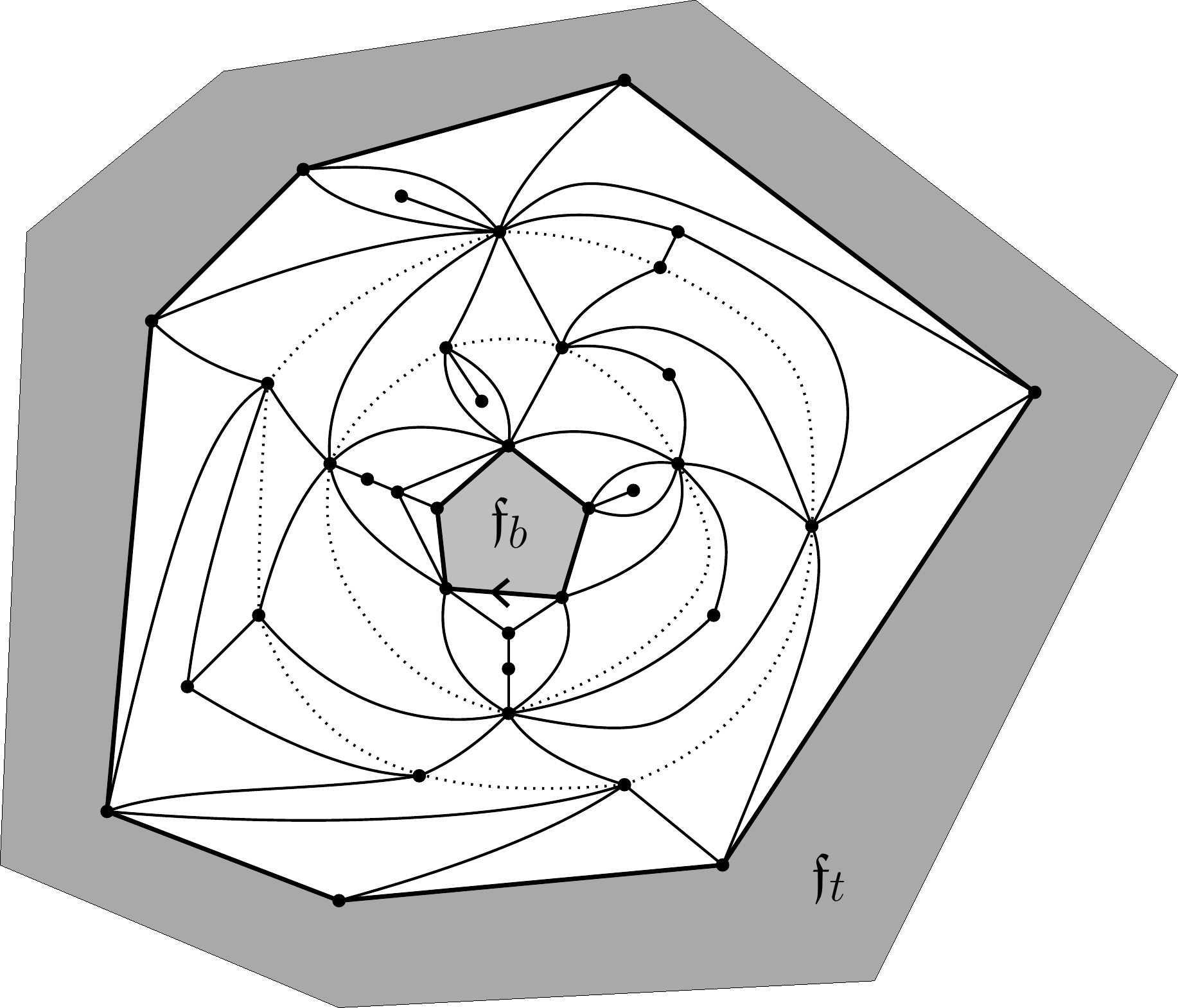}
 \caption{ \label{truncyl} A truncated quadrangulation $\mathcal{Q}$ of the cylinder of height $3$ with boundary sizes $(5,7)$.The two dotted cycles represent $\partial_1\mathcal{Q}$ and $\partial_2\mathcal{Q}$ respectively (see Section \ref{dec-layers} below for the definition of $\partial_k\mathcal{Q}$).}
 \vspace{-5mm}
 \end{center}
 \end{figure}

By definition, the inner faces of $\mathcal{Q}$ are all faces 
except the two distinguished ones. 
The last assertion of Definition \ref{trunc-cyl} shows that the top face and the bottom face do not
play a symmetric role. We will implicitly assume that truncated quadrangulations
of the cylinder of height $h$ are drawn in the plane so that the top face is the unbounded face,
and that they are rooted in such a way that the root edge lies on the bottom cycle and is oriented
clockwise. See Fig.\ref{truncyl} for an example. 

In a way similar to the truncated quadrangulations of Definition \ref{trunc-quad}, 
the triangular face associated with an edge of the bottom cycle
must contain a vertex which does not belong to this cycle. The same holds
for the top cycle --- this is obvious from the last assertion of Definition \ref{trunc-cyl}. 

\subsection{Truncated quadrangulations in the UIPQ}
\label{sec:trunc-UIPQ}

\begin{figure}[!h]
 \begin{center}
 \includegraphics[width=80mm]{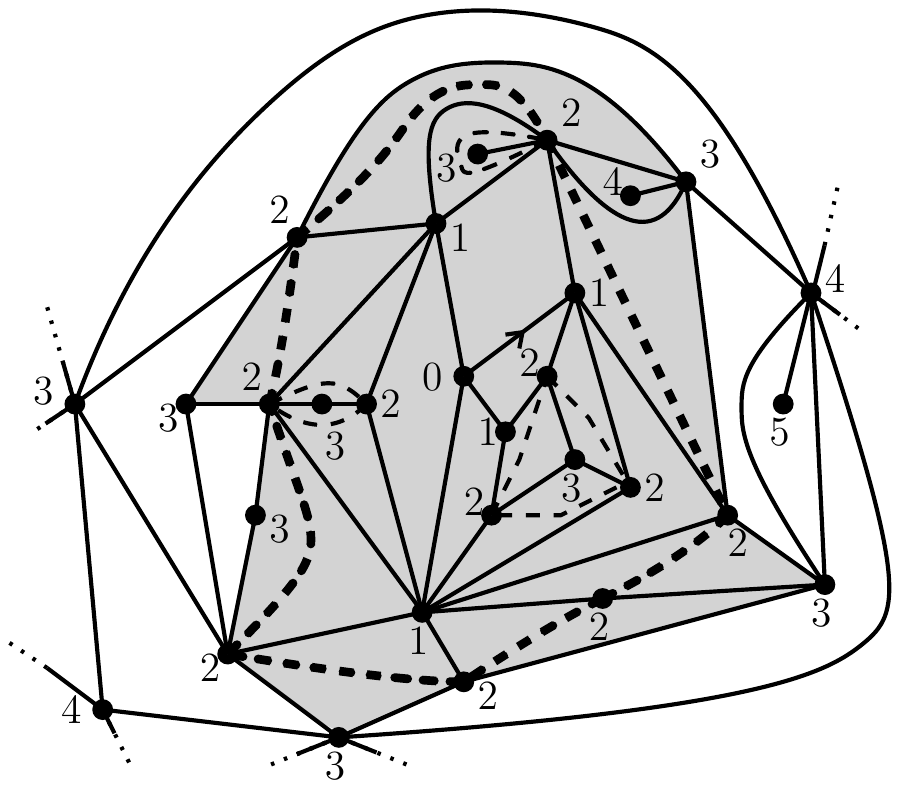}
 \caption{ \label{uipq} The UIPQ near the root vertex. Figures correspond to graph 
 distances from the root vertex. The dashed lines show the cycles made of
 $r$-diagonals, for $r=2$, and the cycle in thick dashed lines is the maximal one. The shaded part is the standard hull of radius $2$. Note that the standard hull
 contains the truncated hull of the same radius, which is the part delimited by the maximal cycle.}
 \vspace{-5mm}
 \end{center}
 \end{figure}

Let us now explain why the definitions of the previous section are relevant to our study
of the UIPQ. 
We label vertices of the UIPQ by their graph distance from the
root vertex. Then the labels of corners incident to a face (enumerated in cyclic order along the boundary of the face) are of the type $k,k-1,k,k-1$
or $k,k+1,k,k-1$ for some integer $k\geq 1$, and the face is called $k$-simple in the second case. Fix an integer $r\geq1$. For every $r$-simple face, we draw 
a diagonal between the two corners labeled $r$ in this face, and these diagonals are called $r$-diagonals. If $v$ is a vertex incident to an $r$-diagonal (equivalently, if $v$ has label $r$ and is incident
to an $r$-simple face), then a simple combinatorial argument shows that the number of $r$-diagonals
incident to $v$ is even --- to be precise, we need to count this number with multiplicities, since
$r$-diagonals may be loops. It follows that
the collection of
all $r$-diagonals can be obtained as the union of a collection of disjoint simple cycles (disjoint here means that 
no edge is shared by two of these cycles). See Fig.\ref{uipq} for an example.

\begin{lemma}
\label{combi-lem}
There is a unique simple cycle made of $r$-diagonals such that the unbounded component of the
complement of this cycle contains no $r$-diagonal and no vertex at distance 
less than or equal to $r$ from the root vertex. This cycle will be called the maximal cycle
made of $r$-diagonals and will be denoted by $\cc_r$.
\end{lemma}

\proof
It suffices to verify that the root vertex lies inside a bounded component of 
the complement of some
cycle made of $r$-diagonals (this cycle may be taken to be simple and then satisfies the properties stated in the lemma).
To this end, consider a geodesic $\gamma$ from the root vertex
to infinity and write $v_r$, resp. $v_{r-1}, v_{r+1}$, for the unique vertex of $\gamma$ at distance 
$r$, resp. $r-1,r+1$, from the root vertex. 
Also write ${v_rv_{r-1}}$, resp. ${v_rv_{r+1}}$, for the edge of $\gamma$
incident to $v_r$ and $v_{r-1}$, resp. to $v_r$ and $v_{r+1}$. Let $k_1$, resp. $k_2$, denote the number of
$r$-diagonals incident to $v_r$ that lie between ${v_rv_{r+1}}$ and ${v_rv_{r-1}}$, 
resp. between ${v_rv_{r-1}}$ and ${v_rv_{r+1}}$, when turning around $v_r$ 
in clockwise order (self-loops are counted twice). An easy combinatorial argument shows that
both $k_1$ and $k_2$ are odd. It follows that there must exist a cycle made 
of $r$-diagonals that starts with an edge lying between ${v_rv_{r+1}}$ and ${v_rv_{r-1}}$ (in clockwise order) and ends with an edge lying between ${v_rv_{r-1}}$ and ${v_rv_{r+1}}$.
Simple topological
considerations now show that the root vertex, and in fact the whole geodesic path $\gamma$
up to vertex $v_{r-1}$ must lie in a bounded component of the complement of this cycle.
\endproof

If we now add all edges of $\cc_r$ to the UIPQ and then remove all edges that lie
in the unbounded component of the complement of $\cc_r$, we get a 
truncated quadrangulation in the sense of Definition \ref{trunc-quad} (with the minor difference that, assuming that 
we keep the same root as in the UIPQ, the root edge
does not belong to the boundary of the distinguished face).
This truncated quadrangulation is called the {\it truncated hull} of radius $r$
and is denoted by $\h^\tr_r$. Its boundary
size (the length of $\cc_r$) is called the perimeter of the hull and denoted by $H_r$. Notice that,
by construction, any vertex belonging to the boundary of the distinguished face is at distance
exactly $r$ from the root vertex. Furthermore, for any vertex $v$ of the UIPQ that does not belong to
$\h^\tr_r$ (equivalently, that lies in the unbounded component of the complement of $\cc_r$) there exists 
a path going from $v$ to infinity that visits only vertices with label at least $r$. This property follows from the fact that any two points of $\cc_r$ are connected by a path that visits only vertices with label at least $r$.

\begin{figure}[!h]
 \begin{center}
 \includegraphics[width=120mm]{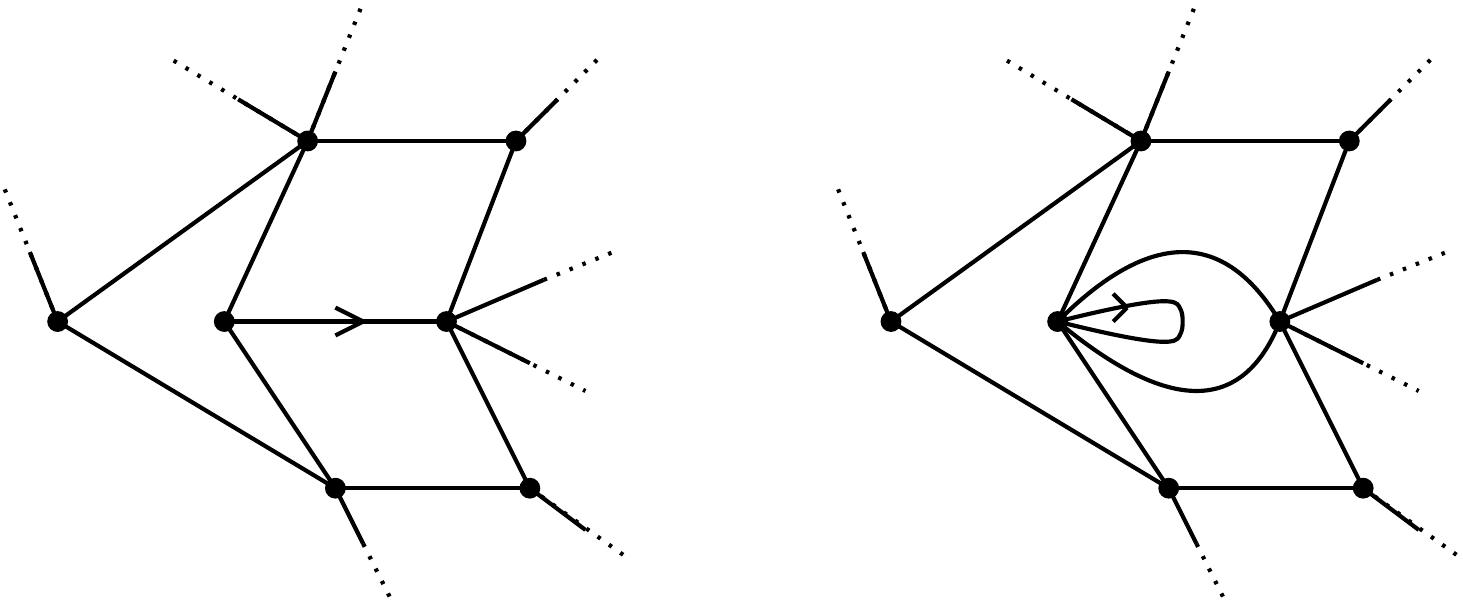}
 \caption{ \label{root-edge} Viewing the hull $\h^\tr_r$ as a quadrangulation of the cylinder:
 The root edge is split in a double edge, and a loop is inserted inside the 
 resulting $2$-gon.}
 \vspace{-5mm}
 \end{center}
 \end{figure}

We may and will sometimes view the truncated hull $\h^\tr_r$ as a quadrangulation of the cylinder:
To this end, we just split the root edge 
into a double edge, and insert a loop (based on the root vertex) inside the
resulting $2$-gon. This yields a truncated quadrangulation of the cylinder of height $r$
with boundary sizes $(1,H_r)$, whose top cycle is $\cc_r$. The root edge is the inserted loop as required 
in our conventions. See Fig.\ref{root-edge} for an illustration. 

Similarly, if $1\leq r<r'$, we can consider the part of the UIPQ that lies between the cycles 
$\cc_r$ and $\cc_{r'}$. More precisely, we add all edges of $\cc_{r}$ and $\cc_{r'}$
to the UIPQ and then remove all edges that lie either inside the cycle $\cc_r$
or outside the cycle $\cc_{r'}$. This gives rise to a quadrangulation of the cylinder of height $r'-r$
whose bottom cycle and top cycle are $\cc_r$ and $\cc_{r'}$ respectively (we in fact need 
to specify the root edge on the bottom cycle, but we will come back to this later). By definition, this
is the {\it annulus} $\cc(r,r')$. We can extend this definition to $r=0$: The annulus 
$\cc(0,r')$ is just the truncated hull $\h^\tr_{r'}$ viewed as a quadrangulation of the cylinder
(we can also say that it is the part of the UIPQ that lies between the
cycles $\cc_0$ and $\cc_{r'}$, if $\cc_0$ consists of the loop added as explained above). 

As an important remark, we note that the truncated hull of radius $r$ is quite different from the (usual) hull of radius $r$ 
considered e.g. in \cite{CLG0,CLG1}, which is 
denoted by $B^\bullet_r(\pp)$ and is obtained by filling in the bounded holes in the ball of radius $r$
(recall that the ball of
radius $r\geq 1$ is obtained as the union of all faces incident to at least one vertex whose
graph distance from the root vertex is at most $r-1$). To avoid any ambiguity, the hull $B^\bullet_r(\pp)$ will
be called the standard hull of radius $r$. The truncated hull can be recovered from the
standard hull by considering the maximal cycle made of $r$-diagonals as explained above.
On the other hand, the standard hull is ``bigger'' than the truncated hull: To recover the 
standard hull from the truncated hull, we need to add the triangles incident to $r$-diagonals
that have been cut when removing the unbounded component of the complement of the maximal cycle,
but also to fill in the bounded holes that may appear when adding these triangles (see Fig.\ref{uipq}
for an example). For future use, we notice that the boundary of the standard hull $B^\bullet_r(\pp)$
is a simple cycle, and that the graph distances of vertices in this cycle to the root vertex
alternate between the values $r$ and $r+1$: Those vertices at graph distance $r$ also
belong to the cycle $\cc_r$, but in general there are other vertices of $\cc_r$ that do not
belong to the boundary of $B^\bullet_r(\pp)$ (see Fig.\ref{uipq}).

\subsection{The skeleton decomposition}
\label{dec-layers}

We will now describe a decomposition of quadrangulations of the cylinder in layers. This is essentially
due to Krikun \cite{Kr} and very similar to the case of triangulations, which is treated in \cite{Kr0,CLG}.
For this reason, we will skip some details.

Let us fix a quadrangulation $\mathcal{Q}$ of the cylinder of height $h\geq 2$ with boundary sizes $(p,q)$.
Assign to each vertex a label equal to its distance from the bottom boundary. Let 
$k\in\{1,2,\ldots,h-1\}$, and consider all diagonals connecting  corners labeled $k$ in 
$k$-simple faces (defined in exactly the same manner as in the previous section for the UIPQ). As in the case of the UIPQ described above, 
these diagonals form a collection of cycles, and there is a maximal cycle
which is simple and has the property that
the unbounded component of the complement of this cycle contains no vertex 
with label less than or equal to $k$.  Define the hull 
$\mathcal{H}_k(\mathcal{Q})$ by first adding to $\mathcal{Q}$
the edges of this maximal cycle and then removing all edges that
lie in the unbounded component of the complement of the maximal cycle. We obtain a quadrangulation of
the cylinder of height $k$ with boundary sizes $(p,q_k)$, where
$q_k$ denotes the size of the maximal cycle. We write $\mathcal{H}^\tr_k(\mathcal{Q})$
for this quadrangulation of the cylinder, and $\partial_k\mathcal{Q}$ for its top cycle ,
so that $q_k=|\partial_k\mathcal{Q}|$. See Fig.\ref{truncyl} for the cycles $\partial_k\mathcal{Q}$
in a particular example. 

Suppose now that we add to $\mathcal{Q}$ all diagonals drawn in the previous procedure, for every
$1\leq k\leq h-1$ (in other words, we add the cycles $\partial_k\mathcal{Q}$ for every
$1\leq k\leq h-1$), and write $\mathcal{Q}^\bullet$
for the resulting planar map (whose faces, except for the two distinguished 
faces of $\mathcal{Q}$, are either quadrangles or triangles).
For every $1\leq k\leq h$, the $k$-th layer of $\mathcal{Q}$ is obtained as the part 
of $\mathcal{Q}^\bullet$ that lies between the cycles $\partial_{k-1}\mathcal{Q}$ and $\partial_k\mathcal{Q}$, where by convention
$\partial_0\mathcal{Q}$ is the bottom cycle of $\mathcal{Q}$ and $\partial_h\mathcal{Q}$ is the top cycle. We can view this
layer as a quadrangulation of the cylinder of height $1$ with boundary sizes
$(|\partial_{k-1}\mathcal{Q}|,|\partial_k\mathcal{Q}|)$ (except that we have not specified the choice 
of the root edge --- we will come back to this later in the case of interest to us).

\begin{figure}[!h]
 \begin{center}
 \includegraphics[width=56mm]{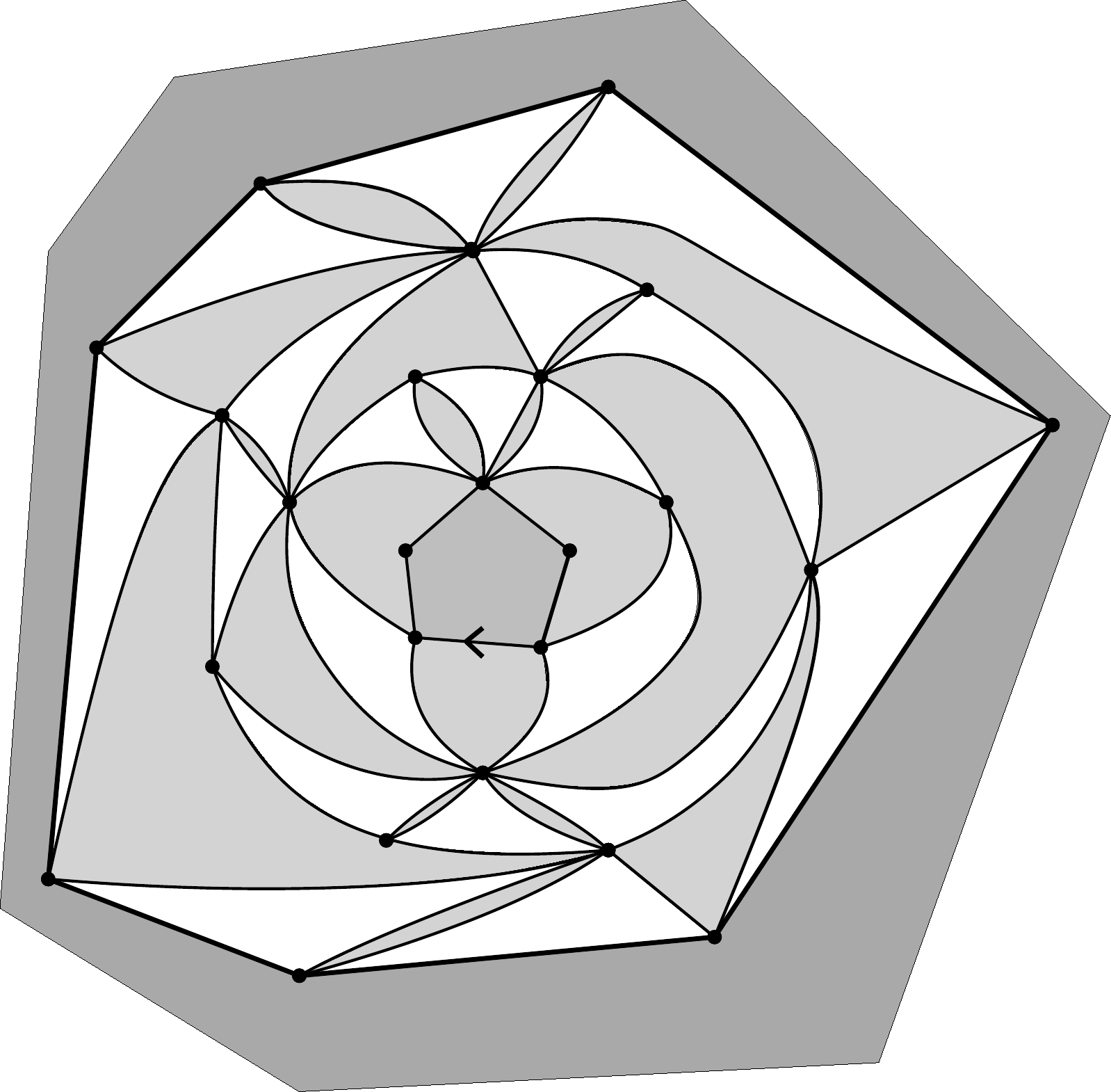}\ \includegraphics[width=62mm]{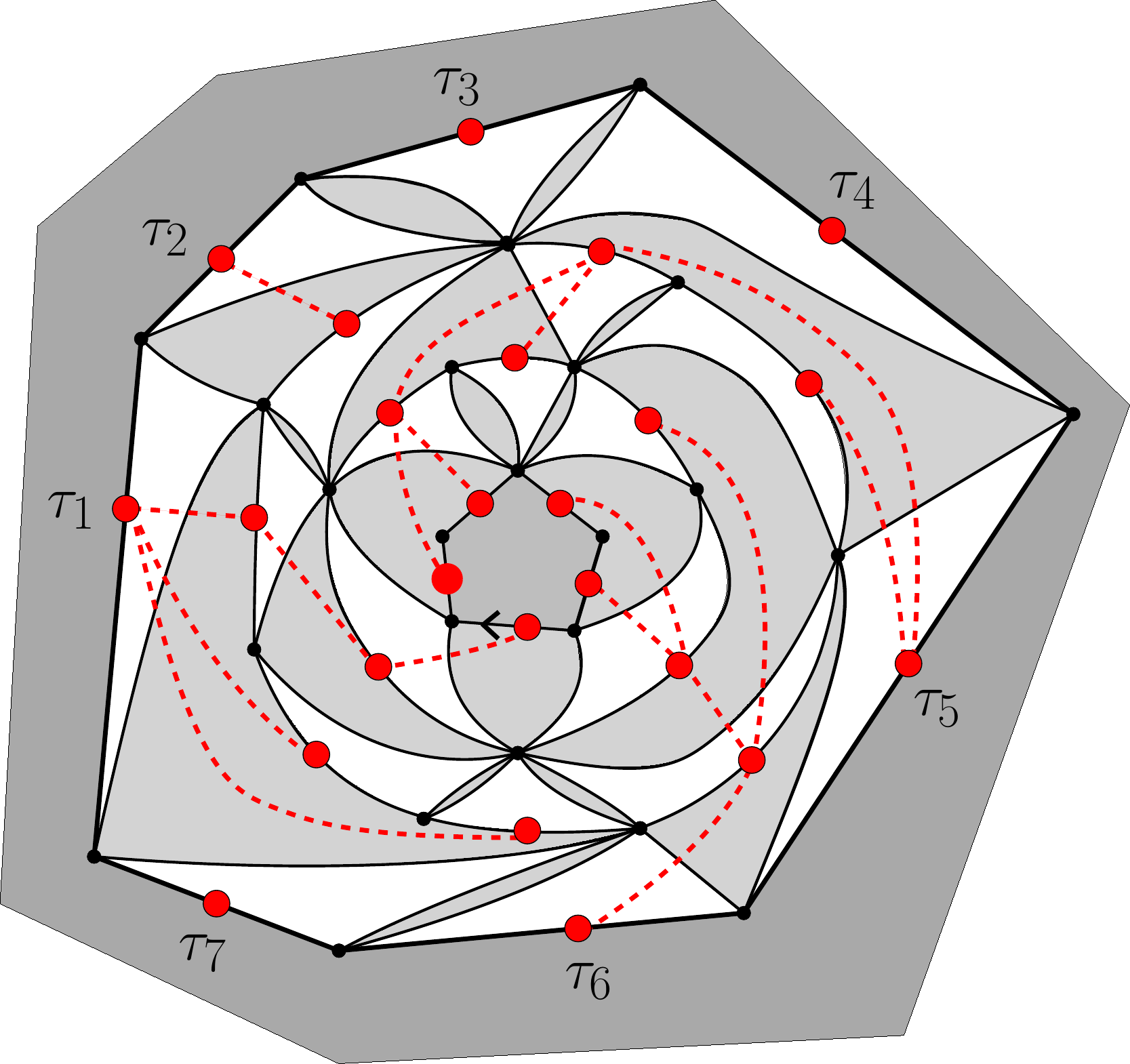}
 \caption{ \label{slots} On the left side, the downward triangles, in white, and the slots, in light grey, in the truncated quadrangulation
 of Fig.\ref{truncyl}. Notice that each edge incident to two downward triangles has been split
 in a double edge, to emphasize the fact that this creates a slot which is a two-gon (whose filling leads to gluing the two sides of the two-gon in Fig.\ref{truncyl}). On the right side, the red dashed lines are the edges of the trees $\tau_1,\ldots,\tau_7$ of the forest 
coding the configuration of downward triangles (notice that $\tau_3,\tau_4,\tau_7$ are trivial trees consisting only of their root vertex). The roots of the trees in the forest are the edges
 of the top cycle, and the trees grow ``toward'' the bottom cycle.}
 \vspace{-5mm}
 \end{center}
 \end{figure}

We will now introduce an unordered forest $\mathcal{F}(\mathcal{Q})$ of (rooted) plane trees that in some sense
describes the configuration of layers. First note that, for every $1\leq k\leq h$, 
each edge of $\partial_k\mathcal{Q}$ is incident to a unique triangle 
of $\mathcal{Q}^\bullet$ whose third vertex lies on $\partial_{k-1}\mathcal{Q}$
(when $k=h$, this is a consequence of the last assertion of Definition \ref{trunc-cyl}, and when 
$k<h$ this follows from the way we constructed the triangles incident to the
top boundary of $\mathcal{H}^\tr_k(\mathcal{Q})$). We call such triangles
{\it downward triangles} of $\mathcal{Q}^\bullet$ (see the left side of Fig.\ref{slots}).
The forest $\mathcal{F}(\mathcal{Q})$ consists of exactly $q$ trees, each tree being associated 
with an edge of $\partial_h\mathcal{Q}$. The vertex set of the forest is the collection of
all edges of $\partial_k\mathcal{Q}$, for $0\leq k\leq h$. The genealogical relation is specified as
follows: The roots of the trees are the edges of $\partial_h\mathcal{Q}$, and, for every $k\in\{0,\ldots,h-1\}$, an edge $e$ of $\partial_{k}\mathcal{Q}$
is a ``child'' of an edge $e'$ of $\partial_{k+1}\mathcal{Q}$ if and only if 
the downward triangle associated with $e'$ (i.e., containing $e'$ in its boundary) is the first one that one encounters 
when turning around $\partial_{k-1}\mathcal{Q}$ in clockwise order, starting from the middle of
the edge $e$. This definition should be clear from
the right side of Fig.\ref{slots}. Notice that edges 
of $\partial_k\mathcal{Q}$ correspond to vertices of the forest $\mathcal{F}(\mathcal{Q})$ at
generation $h-k$, for every $0\leq k\leq h$. The planar structure of each tree in the forest is obviously
induced by the planar structure of $\mathcal{Q}$, see again Fig.\ref{slots}.

We note that the root edge of $\mathcal{Q}$ is a vertex of $\mathcal{F}(\mathcal{Q})$ at generation $h$
and belongs to one of the trees of $\mathcal{F}(\mathcal{Q})$, which we denote by $\tau_1$. We may
then write $\tau_2,\ldots,\tau_q$ for the other trees of of $\mathcal{F}(\mathcal{Q})$ listed in clockwise
order from $\tau_1$. Without risk of confusion, we keep the notation $\mathcal{F}(\mathcal{Q})$
for the ordered forest $(\tau_1,\ldots,\tau_q)$.

\begin{figure}[!h]
 \begin{center}
 \includegraphics[width=12cm]{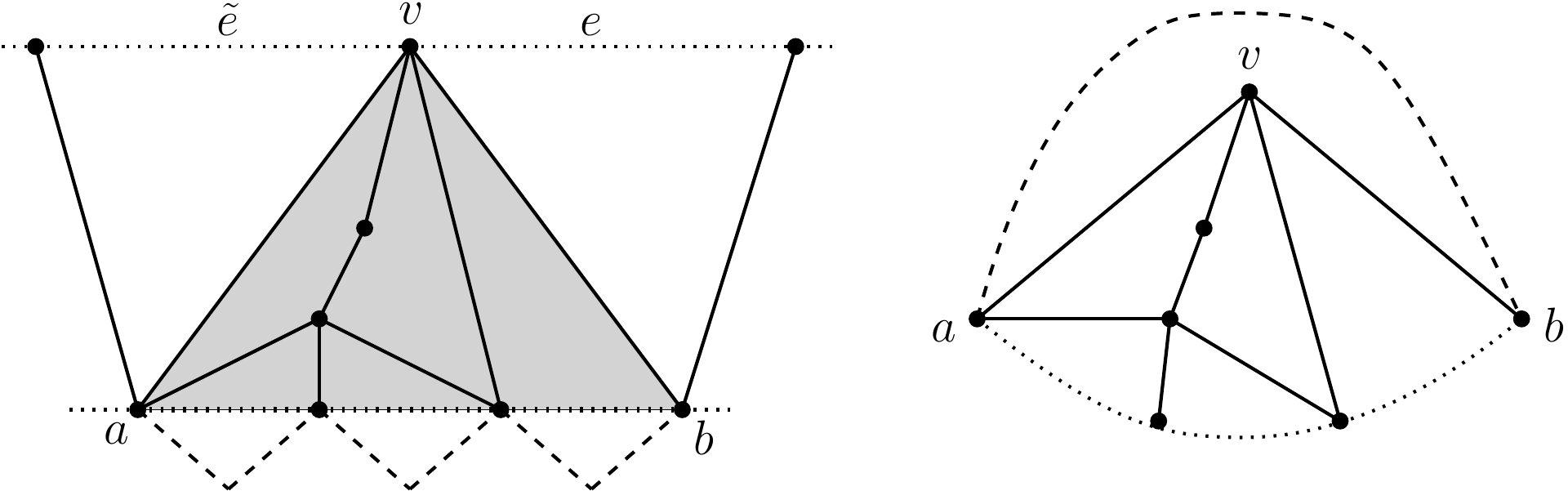}
 \caption{ \label{transform} On the left side, the shaded part corresponds to the slot associated with an edge $e$ of $\partial_k\mathcal{Q}$, such that $c_e=3$. This slot is bounded 
 by the two ``vertical edges'' $av$ and $bv$ (which are incident to the downward triangles associated with $\tilde e$ and $e$ respectively) and by three diagonals (in dotted lines between $a$ and $b$).
 On the right side, this slot is viewed as a truncated triangulation with boundary size $4$ by adding
 the edge between $a$ and $b$ in dashed lines.}
 \vspace{-5mm}
 \end{center}
 \end{figure}

The ordered forest $\mathcal{F}(\mathcal{Q})$ characterizes the combinatorial structure of the
downward triangles in $\mathcal{Q}$. To determine $\mathcal{Q}$ completely, one also needs to specify the
way ``slots'' between two successive downward triangles in a given layer are filled in. 
More precisely, let $e$ be an edge of $\partial_k\mathcal{Q}$, for some
$1\leq k\leq h$, and let $\tilde e$
be the edge of $\partial_k\mathcal{Q}$ preceding $e$ in clockwise order (we discuss below the case
when there is only one edge in $\partial_k\mathcal{Q}$). The part of the $k$-th layer of $\mathcal{Q}$ between the downward triangle 
associated with $\tilde e$ and the downward triangle associated with $e$ produces a slot
with perimeter $c_{e}+2$, where $c_{e}$ is the number of children of $e$
in the forest $\mathcal{F}(\mathcal{Q})$. This slot is said to be associated with $e$ (it is also incident to 
a unique vertex $v$ of $\partial_k\mathcal{Q}$). See the left side of Fig.\ref{transform} for an illustration. If $c_{e}=0$, it may happen that the slot
is empty, if the downward triangles associated with $\tilde e$ and $e$ are adjacent. 
Also notice that when $|\partial_k\mathcal{Q}|=1$, the only edge of $\partial_k\mathcal{Q}$ is a loop,
but there is still an associated slot, which is bounded by the double edge in the boundary
of the downward triangle associated with the unique edge of $\partial_k \mathcal{Q}$, and the edges of $\partial_{k-1}\mathcal{Q}$.

The boundary of the slot associated with $e$   is of the type 
pictured in the left side of Fig.\ref{transform}, where there are $c_{e}$ horizontal edges and the two non-horizontal edges 
are incident to the downward triangles associated with $\tilde e$ and $e$. Strictly speaking, the random
planar map consisting of the part of $\mathcal{Q}^\bullet$ in the slot is not a truncated quadrangulation
with a boundary, but a simple transformation allows us to view it
as a truncated quadrangulation with a boundary of size $c_{\tilde e}+1$: this transformation,
which involves adding
an extra edge, is illustrated in Fig.\ref{transform} (see also Fig.6 in \cite{Kr})  --- to be precise, one should notice that the two vertices $a$ and $b$ in Fig.\ref{transform} 
may be the same if all edges of $\partial_{k-1}\mathcal{Q}$ have the same ``parent'' in $\partial_k \mathcal{Q}$, but our interpretation still goes through. There is therefore a one-to-one 
correspondence between possible fillings of the slot and such truncated quadrangulations. To make
this correspondence precise, we need a convention for the position of the root: we can declare that
in the filling of the slot, the root edge of the truncated quadrangulation corresponds to the 
added extra edge. We notice that, in the special case where $c_{e}=0$, if the truncated
quadrangulation used to fill in the slot is the unique truncated quadrangulation with boundary size $1$
and no quadrangle, this means that the slot is empty so that two sides 
of the downward triangles associated with $\tilde e$ and $e$ are glued together.

Following \cite{CLG}, we say that
a forest $\mathcal{F}$ with a distinguished vertex
is $(h,p,q)$-admissible if
 \begin{enumerate}[(i)]
 \item the forest consists of an ordered sequence $ (\t_{1}, \t_{2}, \ldots , \t_{q})$ 
 of $q$ (rooted) plane trees,
\item the maximal height of these trees is $h$,
\item the total number of vertices of the forest at generation $h$ is  $p$,
\item the distinguished vertex has height $h$ and belongs to $\t_1$.
\end{enumerate}
If  $\mathcal{F}$ is a $(h,p,q)$-admissible forest, we write $\mathcal{F}^*$ for the set all vertices of $\mathcal{F}$ at height strictly less than $h$.
We write $\F^\circ_{h,p,q}$ for the set of
all $(h,p,q)$-admissible forests.

The preceding  discussion yields a bijection between, on the one hand, truncated quadrangulations $\mathcal{Q}$ of the cylinder of height $h$ with boundary sizes $(p,q)$, and, on the other hand, 
pairs consisting of a
$(h,p,q)$-admissible 
forest $\mathcal{F}$ and a collection $(M_v)_{v\in\mathcal{F}^*}$
 such that, for every $v\in\mathcal{F}^*$,  $M_v$ is
a truncated quadrangulation with boundary size $c_v+1$, if $c_v$ stands for the number of children of $v$
in $\mathcal{F}$. We call this bijection the {\it skeleton decomposition} and we say that $\mathcal{F}$ is the {\it skeleton} of the quadrangulation $\mathcal{Q}$. 

It will also be convenient to use the notation $\F_{h,p,q}$ for the set of all 
(ordered) forests, with {\bf no} distinguished vertex, that satisfy properties (i),(ii),(iii) above. 
If $\f\in \F_{h,p,q}$, we keep the notation $\mathcal{F}^*$ for the set all vertices of $\mathcal{F}$ at height strictly less than $h$.

We also set, for every $p\geq 1$, $h\geq 1$,
$$\F^\circ_{h,p}=\bigcup_{q\geq 1} \F^\circ_{h,p,q}\;,\qquad \F_{h,p}=\bigcup_{q\geq 1} \F_{h,p,q}.$$

We conclude this section with a useful observation about connections between the truncated
hull and the standard hull of the UIPQ. Consider two integers $u$ and $r$ with $1\leq u<r$.
Recall that the truncated hull 
$\h^\tr_r$ is viewed as a truncated quadrangulation
of the cylinder of height $r$, whose top cycle is $\cc_r$. Write $\f^\circ_{(r)}$ for the skeleton of this truncated quadrangulation, and also consider the
cycle $\partial_u\h^\tr_r$, which by construction coincides with $\cc_u$. Vertices of this cycle 
are at distance $u$ from the root vertex of the UIPQ, and may or may not belong to
the boundary of the standard hull of radius $u$. However, assuming that $|\cc_u|>1$, if a vertex $v$ of the cycle $\cc_u$
is such that the parents (in the forest
$\f^\circ_{(r)}$) of the two edges of $\cc_u=\partial_u\h^\tr_r$ incident to $v$  are different edges of the
cycle $\cc_{u+1}=\partial_{u+1}\h^\tr_r$, then $v$ must belong to the boundary of the standard hull of radius $u$.
We leave the easy verification of this combinatorial fact to the reader. Notice that this is only
a sufficient condition and that vertices of  $\cc_u$ that do not satisfy this condition may also belong
to the boundary of the standard hull of radius $u$. 

\subsection{Geodesics in the skeleton decomposition}
\label{sec:geodesics}

Consider again a quadrangulation $\mathcal{Q}$ of the cylinder of height $h\geq 2$ with boundary sizes $(p,q)$.
Let $v$ be a vertex of $\partial_k\mathcal{Q}$ for some $k\in\{1,\ldots,h\}$. We assume that $|\partial_k\mathcal{Q}|\geq 2$.
Then $v$ is incident to two downward triangles which both contain an edge of $\partial_k\mathcal{Q}$ and a vertex of $\partial_{k-1}\mathcal{Q}$. Each of these triangles
has an edge incident both to $v$ and to a vertex of $\partial_{k-1}\mathcal{Q}$, and these two edges (which may be the same if the slot incident to $v$
in the $k$-th layer of $\mathcal{Q}$ is empty) are called downward edges
from $v$. If the slot incident to $v$ is nonempty we can in fact define the left downward edge by declaring that it is the first (downward) edge visited when 
exploring the boundary of the slot in clockwise order starting from a point of $\partial_{k-1}\mathcal{Q}$, and the other downward edge is called the right downward edge
(of course if the slot is empty, the left and right downward edges coincide). We leave it to the reader to adapt these definitions in the case
$|\partial_k\mathcal{Q}|=1$ --- in that case the left and right downward edges form a double edge.

\begin{figure}[!h]
 \begin{center}
 \includegraphics[width=12cm]{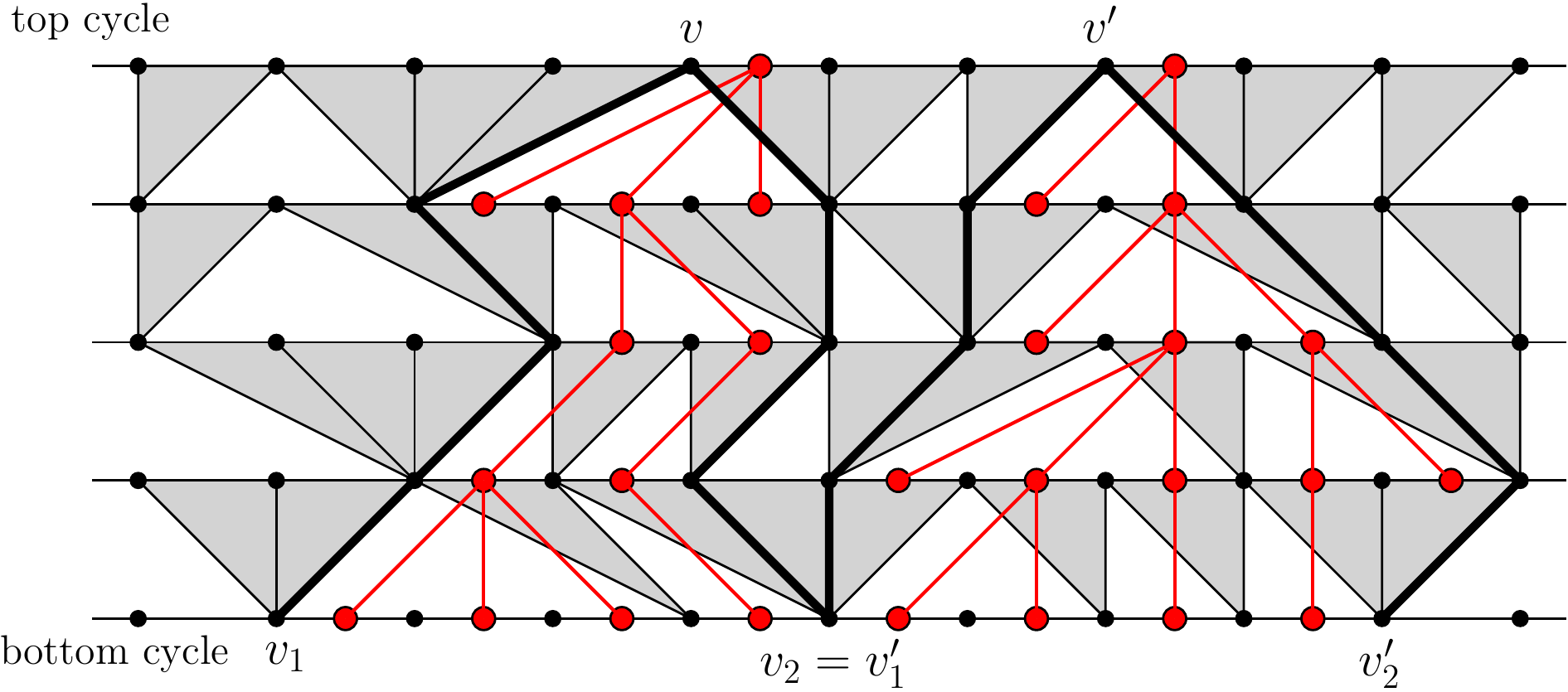}
 \caption{ \label{geode} A portion of a truncated quadrangulation of the cylinder of height $4$. The 
 downward triangles are colored in grey, and the associated slots are in white (they should of
 course be  ``filled in'' by truncated quadrangulations as explained in the text). Two successive trees with maximal height in the
 coding forest are represented in red. The thick black lines are the (left and right) downward geodesics
 from the vertices $v$ and $v'$ of the top cycle associated with the two trees. }
 \vspace{-5mm}
 \end{center}
 \end{figure}

We then define the left downward geodesic from $v$ by saying that we first follow the left downward edge from $v$ to arrive at a vertex
$v'$ of $\partial_{k-1}\mathcal{Q}$, then the left downward edge from $v'$ to a vertex $v''$ of $\partial_{k-2}\mathcal{Q}$, and so on until we reach 
the bottom cycle
$\partial_0\mathcal{Q}$. Similarly we define the right downward geodesic from $v$ by choosing at the first step the right downward edge from $v$, but
then, as previously, following left downward edges from the visited vertices. See Fig.\ref{geode} for an illustration.

Let $N$ be the number of trees 
with maximal height in the skeleton decomposition of $\mathcal{Q}$. Assume that $N\geq 2$, which implies that $|\partial_k\mathcal{Q}|\geq 2$ for every $k\in\{0,1,\ldots,h\}$.
Let $e$ be an edge of $\partial_h\mathcal{Q}$ corresponding to a tree
with maximal height, and let $v$ be the first vertex incident to $e$ in clockwise order around $\partial_h\mathcal{Q}$. Then the left downward geodesic 
(resp. the right downward geodesic) from $v$
hits the bottom cycle at a vertex $v_1$ (resp. at $v_2$) such that the edges of the bottom cycle lying between $v_1$ and $v_2$
in clockwise order are exactly the descendants of $e$ at generation $h$ in the skeleton decomposition. 
See Fig.\ref{geode} for an example.
The concatenation of these two geodesic paths gives a path from $v_1$ to $v_2$ with length $2h$. If we vary the edge 
$e$ among all roots of trees with maximal height, we can concatenate the resulting paths to get a cycle $\cc$ with length $2Nh$,
such that any path from the bottom cycle to the top cycle must visit
a vertex of $\cc$. In particular, if $\mathcal{Q}$ is the annulus $\cc(R,R+h)$ in the UIPQ, with $R\geq 3$, the cycle 
$\cc$ disconnects the ball $B_{R-2}(\pp)$ (or the hull
$B^\bullet_{R-2}(\pp)$) from infinity.

\subsection{Enumeration}
\label{sec:enu}

We rely on the results of Krikun \cite{Kr}. Recall that 
$\Q^{\tr}_{n,p}$ is the set of all truncated quadrangulations with  boundary size $p$ and 
$n$ inner faces (this set is empty if $n<p$). 

\smallskip
Section 2.2 of Krikun \cite{Kr} provides an explicit formula for the generating function
$$U(x,y)=\sum_{p=1}^\infty \sum_{n=1}^\infty \# \Q^{\tr}_{n,p}\ x^ny^p.$$
We will not need this formula, but we record the special case
\begin{equation}
\label{formula-U}
U(\frac{1}{12},y)=\frac{1}{24}\sqrt{(18-y)(2-y)^3}-\frac{1}{2} +\frac{y}{2} -\frac{y^2}{24},
\end{equation}
for $0\leq y<2$. 

As a consequence of the explicit formula for the generating function $U$, we have, for every fixed $p\geq 1$,
\begin{equation}
\label{asym-trunc}
\# \Q^{\tr}_{n,p}\build{\sim}_{n\to\infty}^{} \kappa_p\,n^{-5/2}\,12^n,
\end{equation}
where the constants $\kappa_p$ are determined by the generating function
\begin{equation}
\label{gen-C}
\sum_{p=1}^\infty \kappa_p\,y^p = \frac{128\sqrt{3}}{\sqrt{\pi}}\,\frac{y}{\sqrt{(18-y)(2-y)^3}}.
\end{equation}
for $0\leq y<2$. We again refer to \cite[Section 2.2]{Kr} for these results. 
From \eqref{gen-C} and standard singularity analysis \cite[Corollary VI.1]{FS}, we get
\begin{equation}
\label{asymp-C}
\kappa_p \build{\sim}_{p\to\infty}^{} \frac{64\sqrt{3}}{\pi \sqrt{2}}\, \sqrt{p}\,2^{-p}.
\end{equation}
We also note that
$$\kappa_1=\frac{32}{\sqrt{3\pi}}.$$

\subsection{The distribution of hulls}

Fix integers $n$ and $p$ with $n\geq p$.
Let $\QQ^{(n)}_p$ be uniformly distributed over $\Q^\tr_{n,p}$ and given with a distinguished vertex 
chosen uniformly at random. Let $r\geq 1$. If the height (distance from the boundary) of this
distinguished vertex is at least $r+1$, 
we can make sense
of the hull $\h_r(\QQ^{(n)}_p)$. To this end, we label each vertex by its graph distance from the
boundary of the distinguished face, and we proceed in a way very similar to
the case of the UIPQ discussed in Section \ref{sec:trunc-UIPQ}. We consider all diagonals
connecting corners labeled $r$ in $r$-simple faces (of type $r-1,r,r+1,r$), and the maximal cycle 
made of these diagonals, which has the property that the connected component of the complement
of this cycle containing the distinguished vertex contains only vertices whose label is greater than $r$.
We then add to $\QQ^{(n)}_p$ the edges of this maximal cycle, and remove all edges lying in the
connected component of the complement
of this cycle containing the distinguished vertex. In this way, we obtain 
the hull $\h_r(\QQ^{(n)}_p)$, and it is easy to verify that $\h_r(\QQ^{(n)}_p)$ is a quadrangulation of the cylinder of 
height $r$ (the size of its bottom cycle is $p$). If the height of the distinguished vertex
is smaller than or equal to $r$, the preceding definition no longer makes sense, but by convention
we define $\h_r(\QQ^{(n)}_p)=\dagger$ to be some ``cemetery point'' added 
to the set of all quadrangulations of the cylinder of 
height $r$.

The next lemma, which  is an analog of Lemma 2 in \cite{CLG}, shows that the distribution of $\h_r(\QQ^{(n)}_p)$
has a limit  when $n\to\infty$. We let $\QQ$ be a fixed quadrangulation of the cylinder of height $r$
with boundary sizes $(p,q)$. This quadrangulation is
coded by an $(r,p,q)$-admissible forest $\f$ and a collection $(M_v)_{c\in\f^*}$, such that,
for every $v\in\f^*$, $M_v$ is a truncated quadrangulation with boundary size $c_v+1$. 
Let
$\mathsf{Inn}(M_v)$ denote the number of inner faces of 
$M_v$.

\begin{lemma}
\label{distri-hull}
We have
\begin{equation}
\label{basic}
\lim_{n\to\infty} \P(\h_r(\QQ^{(n)}_p)=\QQ)
= \frac{2^{q}\kappa_q}{2^{p}\kappa_p}\,\prod_{v\in\f^*} \Bigg(\theta(c_v)\,\frac{12^{-\mathsf{Inn}(M_v)}}{Z(c_v+1)}\Bigg)
\end{equation}
where, for every $k\geq1$,
$$Z(k)= \sum_{n=k}^\infty \#\Q^\tr_{n,k}\, 12^{-n},$$
and $\theta$ is the critical offspring distribution defined by
$$\theta(k)=6\cdot 2^{k}\,Z(k+1).$$
The generating function of $\theta$ is given for $0\leq y<1$ by
\begin{equation}
\label{def-gtheta}
g_\theta(y)= 1- \frac{8}{\Big(\sqrt{\frac{9-y}{1-y}}+2\Big)^2 -1}.
\end{equation}
\end{lemma}

We note that the property $Z(k)<\infty$ follows from \eqref{asym-trunc}.

\proof We proceed in a very similar way to the proof of Lemma 2 in \cite{CLG}. Let $N$
be the number of inner faces of $\QQ$, which is also the total number
of vertices of $\QQ$ (by Euler's formula). We observe that the property $\h_r(\QQ^{(n)}_p)=\QQ$
holds if and only if $\QQ^{(n)}_p$ is obtained from $\QQ$ by gluing on the top boundary of $\QQ$
an arbitrary truncated quadrangulation with boundary size $q$ and 
$n-(N-q)$ inner faces (for this gluing to make sense
we need to specify an edge of the top boundary of $\QQ$, which can be the root of the 
first tree in the forest $\f$), and if the distinguished vertex of $\QQ^{(n)}_p$ is chosen among the
inner vertices of this truncated quadrangulation. Noting that $\QQ^{(n)}_p$
has $n+1$ vertices, it follows that
$$\P(\h_r(\QQ^{(n)}_p)=\QQ)
= \frac{\#\Q^\tr_{n-(N-q),q}}{\#\Q^\tr_{n,p}} \times \frac{n+1-N}{n+1}.$$
Using \eqref{asym-trunc}, we get
\begin{equation}
\label{basic1}
\lim_{n\to\infty} \P(\h_r(\QQ^{(n)}_p)=\QQ) = \frac{\kappa_q}{\kappa_p} \,12^{-N+q}.
\end{equation}

Simple combinatorics shows that the number of inner faces of $\QQ$ can be written as 
\begin{equation}
\label{formula-inner}
N= p+ \sum_{v \in \mathcal{F}^*}  (\mathsf{Inn}( M_v) - c_v).
\end{equation}
So the right-hand side of \eqref{basic1} is also equal to
$$\frac{\kappa_q}{\kappa_p}\,12^{q-p} \prod_{v\in \f^*} (12^{c_v}\,12^{- \mathsf{Inn}( M_v)}).$$
It is now straightforward to verify that the last quantity is equal to the
right-hand side of \eqref{basic}. Just observe that
$$12^{c_v}\,Z(c_v+1) = 6^{c_v-1}\,\theta(c_v)$$
and notice that
$$\sum_{v\in \f^*} (c_v-1)= p-q.$$

To see that $\theta$ is an offspring distribution, we rely on \eqref{formula-U}, which shows that
the generating function of $\theta$ is 
\begin{align*}
g_\theta(y):=\sum_{k=0}^\infty \theta(k)\,y^k&=6\sum_{k=0}^\infty\sum_{n=k+1}^\infty y^k 2^k 12^{-n}\# \Q^\tr_{n,k+1}\\
&=\frac{6}{2y}\,U(\frac{1}{12},2y)\\
&=\frac{1}{2y} \Big(\sqrt{(9-y)(1-y)^3} - 3 +6y - y^2\Big),
\end{align*}
in agreement with Theorem 2 of \cite{Kr}. 
Since $g_\theta(1)=1$, $\theta$
is a probability distribution, and
the fact that $\theta$ is critical is obtained by checking that $g_\theta'(1)=1$.

Finally, a somewhat tedious calculation shows that the formula for $g_\theta$ in the last
display is equivalent to the one given in the statement of the lemma. The latter is more
convenient to compute iterates of $g_\theta$, as we will see below in formula \eqref{iterate-g}. 
\endproof

From the explicit form of $g_\theta$, we have
$$g_\theta(1-x)=1-x+\sqrt{2}\,x^{3/2}+ O(x^2)$$
as $x\downarrow 0$. By singularity analysis, it follows that
\begin{equation}
\label{equiv-theta}
\theta(k)\build{\sim}_{k\to\infty}^{} \frac{3\sqrt{2}}{4\sqrt{\pi}}\,k^{-5/2}.
\end{equation}

\rem The offspring distribution $\theta$ appears in the seemingly different context
of labeled trees. Consider a critical Galton-Watson tree with geometric offspring distribution
with parameter $1/2$. Given the tree, assign labels to vertices by declaring 
that the label of the root is $0$ and that label increments on different edges are 
independent and uniformly distributed over $\{-1,0,1\}$. Let $N$ be the 
number of vertices labeled $-1$ whose (strict) ancestors all have nonnegative labels.
Then $N$ is distributed according to $\theta$ (see \cite[Proof of Theorem 5.2]{CMM}).
Via Schaeffer's bijection relating plane quadrangulations to labeled trees, this 
interpretation of $\theta$ is in fact closely related to Lemma \ref{distri-hull}.

\medskip

We define, for every $p\geq 1$,
\begin{equation}
\label{def-h}
h(p):= \frac{1}{p}\,2^{p}\,\kappa_p.
\end{equation}

\begin{lemma}
\label{proba-meas}
Let $p\geq 1$. The formula
$$\mu_{r,p}(\f):=  \frac{h(q)}{h(p)}\,\prod_{v\in\f^*} \theta(c_v)\;,\ \f\in \F_{r,p,q},\,q\geq 1$$
defines a probability measure on $\F_{r,p}$. Consequently, the formula
$$\mu^\circ_{r,p}(\f):=  \frac{2^{q}\kappa_q}{2^{p}\kappa_p}\,\prod_{v\in\f^*} \theta(c_v)\;,\ \f\in \F^\circ_{r,p,q},
\,q\geq 1$$
defines a probability measure on $\F^\circ_{r,p}$. 
\end{lemma}

\proof
Let $\Pi$ be the generating function of the sequence $(h(k))_{k\geq 1}$,
$$\Pi(x):=\sum_{k=1}^\infty h(k)\,x^k.$$
To verify that $\mu_{r,p}$ defines a probability
distribution on the set $\F_{r,p}$, it is enough to
check that $(h(k))_{k\geq 1}$ is an (infinite) stationary measure for the
branching process with offspring distribution $\theta$, or
equivalently that, for every $0<y<1$,
\begin{equation}
\label{stat-meas}
\Pi(g_\theta(y)) - \Pi(g_\theta(0))=\Pi(y).
\end{equation}
From \eqref{gen-C}, we get by integration that
$$\sum_{p=1}^\infty \kappa_p\,\frac{x^p}{p}= \frac{48}{\sqrt{\pi}}\,\Bigg( \sqrt{
\frac{18-x}{3(2-x)}} - \sqrt{3}\Bigg),$$
and, for $0<x<1$,
$$\Pi(x)= \frac{48}{\sqrt{3\pi}}\, \Bigg(\sqrt{\frac{9-x}{1-x}} - 3\Bigg).$$
From this explicit formula and \eqref{def-gtheta}, the desired identity 
\eqref{stat-meas} follows at once. 

Once we know that $\mu_{r,p}$ is a probability distribution on
$\F_{r,p}$, the fact that $\mu^\circ_{r,p}$ is a probability distribution on $\F^\circ_{r,p}$
follows easily. First note that
$$\tilde \mu((\f,v)):= \frac{1}{p}\,\mu_{r,p}(\f)$$
defines a probability distribution on the set of all pairs $(\f,v)$
consisting of a forest $\f\in\F_{r,p}$ and a distinguished vertex
$v$ of $\f$ at generation $r$. Then notice that $\mu^\circ_{r,p}$ is just the
push forward of $\tilde\mu$ under the mapping 
$(\f,v)\mapsto \f'$, where $\f'$ is obtained by circularly permuting the trees
of $\f$ so that $v$ belongs to the first tree of the forest. This completes the 
proof. \endproof

Let us consider now the UIPQ. Recall our notation $\h^\tr_r$ for the
truncated hull of radius $r$, and $H_r$ for the perimeter of $\h^\tr_r$, which is also the 
length of the cycle $\cc_r$. 
For every $r,p\geq 1$, let $\C_{r,p}$ be the set of all 
truncated quadrangulations of the cylinder of height $r$ with bottom boundary size $p$
and arbitrary top boundary size. If $\QQ\in \C_{r,p}$ and the size of the top
boundary of $\QQ$ is $q$, we set
$$\Delta_{r,p}(\QQ)= \frac{2^{q}\kappa_q}{2^{p}\kappa_p}\,\prod_{v\in\f^*} \Bigg(\theta(c_v)\,\frac{12^{-\mathsf{Inn}(M_v)}}{Z(c_v+1)}\Bigg),$$
where 
$(\f,(M_v)_{v\in\f^*})$ is the skeleton decomposition of $\QQ$. 

\begin{corollary}
\label{law-hull}
$\Delta_{r,p}$ is a probability measure on $\C_{r,p}$. Furthermore, the 
distribution of $\h^\tr_r$ is $\Delta_{r,1}$.
\end{corollary}

\proof The fact that $\Delta_{r,p}$ is a probability measure on $\C_{r,p}$ 
readily follows from the second assertion of Lemma \ref{proba-meas}, noting that, by the
very definition of $Z(k)$, we have
$$\sum_{n=1}^\infty\#\Q^\tr_{n,k+1}\, \frac{12^{-n}}{Z(k+1)}= 1.$$
To get the second assertion of the corollary, let $\Q_n$ stand for the set of
all (rooted) planar quadrangulations with $n$ faces. Via the transformation 
that consists in splitting the root edge to get a double edge, and then inserting
a loop inside the resulting $2$-gon (as in Fig.\ref{root-edge}), the set $\Q_{n}$ is canonically identified 
to $\Q^\tr_{n+1,1}$. From the local convergence of planar quadrangulations
to the UIPQ \cite{Kr}, we deduce that the distribution of the hull of radius 
$r$ in a uniformly distributed quadrangulation in $\Q^\tr_{n,1}$ (equipped with a 
distinguished uniformly distributed vertex) converges to the
distribution of the hull
of radius $r$ in the UIPQ. The second assertion of the corollary now follows from
Lemma \ref{distri-hull}. 
\endproof

\begin{corollary}
\label{hull-peri}
The distribution of $H_r$
is given by
$$\P(H_r=p)=\frac{h(p)}{h(1)}\,\P_p(Y_r=1),\quad p\geq 1,$$
where $(Y_n)_{n\geq 0}$ denotes a Galton-Watson branching
process with offspring distribution $\theta$ that starts from
$p$ under the probability measure $\P_p$.
\end{corollary}

\proof
Let $\f_{(r)}^\circ$ be the skeleton of the hull $\h^\tr_r$
viewed as a quadrangulation of the cylinder of height $r$. As a direct consequence of the second assertion of Corollary \ref{law-hull}, 
$\f_{(r)}^\circ$ is distributed as $\mu^\circ_{r,1}$. Define $\f_{(r)}$ from $\f^\circ_{(r)}$ by
``forgetting'' the distinguished vertex and  
applying a uniform random circular permutation to the trees in the sequence. Arguing as 
in the end of the proof of Lemma \ref{proba-meas}, it follows that 
$\f_{(r)}$ is distributed according to $\mu_{r,1}$.

Since $H_r$ is just the number of trees in the forest $\f_{(r)}$, we have
$$\P(H_r=p)=\sum_{\f\in \F_{r,1,p}} \P(\f_{(r)}=\f)
= \frac{h(p)}{h(1)}\, \sum_{\f\in \F_{r,1,p}} \prod_{v\in\f^*} \theta(c_v),$$
and the desired result follows. 
\endproof

\noindent{\bf Remark.} One can interpret the distribution of $H_r$ as the limit when $T\to\infty$
of the distribution at time $T-r$ of a Galton-Watson process $X$ with offspring distribution $\theta$ started 
from $X_0=1$ and conditioned on the event $\{X_T=1\}$. This suggests that one may code the combinatorial
structure of downward triangles in the whole UIPQ (and not only in a hull of fixed radius) by an infinite 
tree, which could be viewed as the genealogical tree for a Galton-Watson process with offspring distribution $\theta$, indexed 
by nonpositive integer times and conditioned to be equal to $1$ at time $0$. This interpretation will not
be needed in the present work and we omit the details. 

\medskip

Let us now fix integers $0\leq u<r$. As explained earlier, the annulus $\cc({u,r})$ is the part of the 
UIPQ that lies between the cycles $\cc_u$ and $\cc_r$ --- recall our convention for
$\cc_0$ from Section \ref{sec:trunc-UIPQ} ---
and $\cc({u,r})$ is viewed as a truncated quadrangulation of the cylinder
of height $r-u$ with boundary sizes $(H_u,H_r)$. We now specify the
root edge of $\cc({u,r})$, by declaring that it corresponds to the root of the tree,
in the skeleton decomposition of $\h^\tr_u$, that
carries the root edge of the UIPQ (of course when $u=0$, the root edge is the unique edge of $\cc_0$).

 Let $\f^\circ_{u,r}$ be the skeleton of $\cc(u,r)$, which is a random element 
 of $\cup_{p\geq 1}\F^\circ_{r-u,p}$. 
It will be convenient to introduce also the forest $\f_{u,r}$ (in $\cup_{p\geq 1}\F_{r-u,p}$)
obtained from $\f^\circ_{u,r}$ by first ``forgetting'' the distinguished vertex and then 
applying a uniform random circular permutation to the trees in the sequence.

\begin{corollary}
\label{law-hull-annulus}
Let $p\geq 1$. The conditional distribution of $\f_{u,r}$ knowing that
$H_u=p$ is $\mu_{r-u,p}$. 
\end{corollary}

\proof 
Recall the notation $\f^\circ_{(r)},\f_{(r)}$ introduced in the previous proof.
We notice that, if $(\tau_1,\ldots,\tau_{H_r})$ are the trees in the forest $\f^\circ_{(r)}$,
the trees in the forest $\f^\circ_{u,r}$ are just $(\tau_1^{[r-u]},\ldots,\tau_{H_r}^{[r-u]})$,
where the notation $\tau_i^{[r-u]}$ refers to the tree $\tau_i$ truncated at generation 
$r-u$. 
It follows that $\f_{u,r}$ can be assumed to be equal to the forest 
$\f_{(r)}$ truncated at generation $r-u$. Note that $H_u$ is just the number 
of vertices of $\f_{(r)}$ at generation $r-u$.

Let $q\geq 1$ and $\g\in\F_{r-u,p,q}$. We have
$$\P(\f_{u,r}=\g)= \sum_{\f\in \F_{r,1} : \f^{[r-u]}=\g} \P(\f_{(r)}=\f),$$
using the notation $\f^{[r-u]}$ for the forest $\f$ truncated at
generation $r-u$. It follows that
\begin{align*}
\P(\f_{u,r}=\g)&= \frac{h(q)}{h(1)} \sum_{\f\in \F_{r,1} : \f^{[r-u]}=\g} \prod_{v\in\f^*} \theta(c_v)\\
&= \frac{h(q)}{h(1)}\,\prod_{v\in\g^*} \theta(c_v)\,\sum_{\tilde\f\in\F_{u,1,p}} \prod_{v\in\tilde\f^*}\theta(c_v),
\end{align*}
where we just use the fact that a forest $\f\in \F_{r,1}$ such that $\f^{[r-u]}=\g$ is obtained by
``gluing'' a forest of $\F_{u,1,p}$ to the $p$ vertices of $\g$ at generation $r-u$. As 
in Corollary \ref{hull-peri} and its proof, we have
$$\sum_{\tilde\f\in\F_{u,1,p}} \prod_{v\in\tilde\f^*}\theta(c_v) = \P_p(Y_u=1)= \frac{h(1)}{h(p)}\,\P(H_u=p).$$
So we get
$$\P(\f_{u,r}=\g)= \frac{h(q)}{h(p)} \prod_{v\in\g^*} \theta(c_v) \times \P(H_u=p)
= \P(H_u=p)\,\mu_{r-u}(\g).$$
This completes the proof. \endproof

\subsection{The law of the perimeter of hulls}

We give a more explicit formula for the distribution of $H_r$.

\begin{proposition}
\label{hull-peri2}
We have, for every $r\geq 1$ and $p\geq1$,
\begin{equation}
\label{law-peri}
\P(H_r=p)=K_r\,\kappa_p\,(2\pi_r)^p
\end{equation}
where
\begin{align*}
\pi_r&=1- \frac{8}{(3+2r)^2 -1}=\frac{r(r+3)}{(r+1)(r+2)}\\
K_r&= \frac{32}{3\kappa_1}\,\frac{3+2r}{((3+2r)^2-1)^2}\;\frac{1}{\pi_r}
=\frac{2}{3\kappa_1}\,\frac{2r+3}{r(r+1)(r+2)(r+3)}.
\end{align*}
Consequently, there exist positive constants $M_1,M_2$ and $\rho$ such that, for every $a>0$,
for every integer $r\geq 1$,
\begin{equation}
\label{bound-H1}
\P(H_r\geq a\,r^2)\leq M_1\,e^{-\rho a}
\end{equation}
and 
\begin{equation}
\label{bound-H2}
\P(H_r\leq a\,r^2)\leq M_2\,a^{3/2}.
\end{equation}
\end{proposition}

We notice that $K_r\sim (4/(3\kappa_1))\,r^{-3}$ as $r\to\infty$ and recall that $\kappa_1=\frac{32}{\sqrt{3\pi}}$. 

\proof
We rely on the formula of Corollary \ref{hull-peri}. Recalling that $(Y_n)_{n\geq 0}$ denotes a Galton-Watson branching
process with offspring distribution $\theta$ that starts from
$p$ under the probability measure $\P_p$, and using formula \eqref{def-gtheta}, we obtain that
the generating function of $Y_r$ under $\P_1$ is
\begin{equation}
\label{iterate-g}
g^{(r)}_\theta(y)= \underbrace{g_\theta\circ\cdots\circ g_\theta}
_{r\;{\rm times}}(y)=  1- \frac{8}{\Big(\sqrt{\frac{9-y}{1-y}}+2r\Big)^2 -1}.
\end{equation}
It follows that
\begin{equation}
\label{extinct-proba}
\P_1(Y_r=0)=g_\theta^{(r)}(0)= 1- \frac{8}{(3+2r)^2 -1}=\pi_r,
\end{equation}
and
\begin{align}
\label{law-Y}
\P_p(Y_r=1)&=\lim_{x\da 0} \frac{1}{x} \Big(\E_p[x^{Y_r}] - \P_p(Y_r=0)\Big)\nonumber\\
&=\lim_{x\da 0} \frac{1}{x} \Big(g^{(r)}_\theta(x)^p- g^{(r)}_\theta(0)^p\Big)\nonumber\\
&=p\, g^{(r)}_\theta(0)^{p-1} \times \frac{64}{3}\,\frac{3+2r}{((3+2r)^2-1)^2}\nonumber\\
&=\frac{64}{3}\,p\,\frac{3+2r}{((3+2r)^2-1)^2}\,\Bigg(1- \frac{8}{(3+2r)^2 -1}\Bigg)^{p-1}.
\end{align}
Since
$$\frac{h(p)}{h(1)}= \frac{1}{p}\,\frac{2^p\kappa_p}{2\,\kappa_1},$$
Corollary \ref{hull-peri} and \eqref{law-Y} lead to formula
\eqref{law-peri}. Finally, the bounds \eqref{bound-H1} and \eqref{bound-H2} are
simple consequences of this explicit formula and the asymptotics \eqref{asymp-C}
for the constants $\kappa_p$. To derive \eqref{bound-H1}, we observe that
we can find a constant $\eta>0$ such that
$\P(H_r\geq a\,r^2)$ is bounded above by a constant times
$$r^{-3}\sum_{p>ar^2} \sqrt{p}\times e^{-\eta p/r^2} \leq \hbox{Cst.}\,r^{-3}\int_{ar^2/2}^\infty \sqrt{x}e^{-\eta x/r^2} \mathrm{d}x
=\hbox{Cst.}\times \int_{a/2}^\infty \sqrt{y}e^{-\eta y} \mathrm{d}y,$$
and the proof of \eqref{bound-H2} is even easier just bounding $\pi_r$ by $1$. 
\endproof

\subsection{A conditional limit for branching processes}

We keep the notation $(Y_n)_{n\geq 0}$ for
a branching process with offspring distribution 
$\theta$, which starts at $p$ under the probability measure $\P_p$.

\begin{lemma}
\label{limit-branching}
We have
$$\P_1(Y_r\not =0)\build{\sim}_{r\to\infty}^{} \frac{2}{r^2}$$
and the distribution of $r^{-2}Y_r$ under
$\P_1(\cdot\mid Y_r\not =0)$ converges to the distribution with Laplace transform
$$1- \Big(1 +\sqrt{\frac{2}{\lambda}}\Big)^{-2}.$$
\end{lemma}

\proof The first assertion is immediate from  \eqref{extinct-proba}. 
Next, from \eqref{iterate-g} and \eqref{extinct-proba}, we have
$$\E_1\Big[ e^{-\lambda r^{-2}Y_r}\,\mathbf{1}_{\{Y_r\not = 0\}}\Big]
= \frac{8}{(3+2r)^2 -1} - \frac{8}{\Big(\sqrt{1+\frac{8}{1-e^{-\lambda r^{-2}}}}+2r\Big)^2 -1}$$
and it easily follows that
$$\frac{r^2}{2}\,\E_1\Big[ e^{-\lambda r^{-2}Y_r}\,\mathbf{1}_{\{Y_r\not = 0\}}\Big]
\build{\longrightarrow}_{r\to\infty}^{} 1- \Big(1 +\sqrt{\frac{2}{\lambda}}\Big)^{-2},$$
giving the desired result. \endproof

\subsection{An estimate on discrete bridges}
\label{bridge}

In this short section, which is independent of the previous ones, we state an estimate for discrete 
bridges, which plays an important role in the proof of Proposition \ref{main} in the next section.

Let $K\geq 1$ be an integer, and let $(b(0),b(1),\ldots,b(2K))$ be a discrete bridge of length $2K$. This 
means that $(b(0),b(1),\ldots,b(2K))$ is uniformly distributed over sequences $(x_0,x_1,\ldots,x_{2K})$
such that $x_0=x_{2K}=0$ and $|x_i-x_{i-1}|=1$ for every $i=1,\ldots,2K$. It will be convenient 
to define intervals on $\{0,1,\ldots,2K-1\}$ in a cyclic manner: If $i,j\in\{0,1,\ldots, 2K-1\}$,
$[i,j]=\{i,i+1,\ldots,j\}$ as usual if $i\leq j$, but 
$[i,j]=\{i,i+1,\ldots,2K-1,0,1,\ldots,j\}$ if $i>j$. 

Let $c>0$ be a fixed constant and let $r\geq 1$ be an integer. For every integer $k\geq 2$, we let $\pp_{k,K}(r)$
stand for the event where there exist integers $0\leq m_1<m_2<\cdots<m_k<2K$, such that 
$m_i-m_{i-1}\geq cr^2$ for every $2\leq i\leq k$, and $m_1+2K-m_k\geq cr^2$, and,
for every $i,j\in\{1,\ldots,k\}$,
$$b(m_i)+b(m_j) - 2 \,\max\Big( \min_{\ell\in[m_i,m_j]} b(\ell),\min_{\ell\in[m_j,m_i]} b(\ell)\Big) \leq 5\,r.$$

\begin{lemma}
\label{tec-lem}
There exist constants $C>0$ and $\gamma\in(0,1)$, which only depend on $c$, 
such that, for every $r\geq 1$ and every $k\geq 2$,
$$\P(\pp_{k,K}(r))\leq \frac{C}{k}\,\Big(\frac{K}{r^2}\Big)^2\,\gamma^k.$$
\end{lemma}

We postpone the proof to the Appendix.

\section{Lower bound on the size of the separating cycle}
\label{sec:main}

In this section, we prove Proposition \ref{main}, and then explain how
part (i) of Theorem \ref{kri-conj} follows from this result.

\medskip
\noindent{\it Proof of Proposition \ref{main}.}
As a preliminary observation, we note that it is enough to prove that
the stated bound holds for $n$ large enough (and for every $r\geq 1$).
Let $u$ and $w$ be two integers with $0\leq u<w$. Recall the notation
$\f_{u,w}$ for the forest obtained from the skeleton of  $\cc(u,w)$ 
by forgetting the distinguished vertex and then 
applying a uniform random circular permutation to the trees in the forest. By Corollary \ref{law-hull-annulus}, we have
for every $p,q\geq1$ and $\f\in \F_{w-u,p,q}$,
$$\P(\f_{u,w}=\f\mid H_u=p) = \frac{h(q)}{h(p)}\,\prod_{v\in \f^*} \theta(c_v),$$
where $h$ is defined in \eqref{def-h}. By Corollary \ref{hull-peri}, we have
$$\frac{h(q)}{h(p)}=\frac{\P(H_w=q)}{\P(H_u=p)}\,\frac{\P_p(Y_u=1)}{\P_q(Y_w=1)}$$
where we must take $p=1$ if $u=0$.
It follows that
$$\P(\f_{u,w}=\f)= \P(H_w=q)\,\frac{\P_p(Y_u=1)}{\P_q(Y_w=1)}\,\prod_{v\in \f^*} \theta(c_v)
$$
and therefore
\begin{equation}
\label{law-forest}
\P(\f_{u,w}=\f\mid H_w=q)=\frac{\P_p(Y_u=1)}{\P_q(Y_w=1)}\,\prod_{v\in \f^*} \theta(c_v)
=\frac{\varphi_u(p)}{\varphi_w(q)}\,\prod_{v\in \f^*} \theta(c_v),
\end{equation}
with the notation
\begin{equation}
\label{formula-phi}
\varphi_u(p)=\P_p(Y_u=1)=\frac{64}{3}\,p\,\frac{3+2u}{((3+2u)^2-1)^2}\,\pi_u^{p-1},
\end{equation}
by \eqref{law-Y}.
We will apply formula \eqref{law-forest} with $u=nr$ and $w=(n+2)r$ for integers $n,r\geq 1$. 

Let us fix $\alpha>0,\beta>0$ with $\alpha<2<\beta$ and $\beta-2<\alpha$.
Let $c_0>0$ be a constant whose value will be specified later.
Say that a plane tree satisfies property $(P)_r$ if it has at least $c_0r^2$ vertices of generation $2r-1$ 
that have at least one child at generation $2r$.
Thanks to Lemma \ref{limit-branching}, we can choose the constant $c_0>0$ 
small enough
so that, for every $r$ large enough,
the probability for a Galton-Watson tree with offspring distribution $\theta$
to satisfy property $(P)_r$  is greater than
$a_0r^{-2}$, for some other constant $a_0>0$.
Let $\delta$ be another constant, with $0<\delta<\alpha$.  We write
$\Theta_{n,r}$ for the collection of all forests in
$$\bigcup_{p,q\geq 1} \F_{2r,p,q}$$
having at least $n^\delta$ trees that satisfy property $(P)_r$. Our first goal is to find
an upper bound for 
$$\P(\f_{nr,(n+2)r}\notin \Theta_{n,r}).$$

 From \eqref{bound-H1}, we have
 \begin{equation}
\label{first-esti}
\sup_{r\geq 1}\Bigg(\max\Big(\P(H_{nr}>n^\beta r^2),\P(H_{(n+2)r} > n^\beta r^2)\Big)\Bigg) = O(e^{-n^\ve})
\end{equation}
as $n\to\infty$, with a constant $\ve>0$ that depends only on $\beta$. On the other hand, by \eqref{bound-H2}, we have also
\begin{equation}
\label{2-esti}
\P(H_{(n+2)r}<n^\alpha r^2)\leq M_2\,n^{\frac{3}{2}(\alpha-2)},
\end{equation}
and the constant $M_2$ does not depend on $r$. 

We then restrict our attention to
\begin{align}
\label{decompo}
&\P(\f_{nr,(n+2)r}\notin \Theta_{n,r},\,n^\alpha r^2\leq H_{(n+2)r}\leq n^\beta r^2,\,H_{nr}\leq n^\beta r^2)\nonumber\\
&= \sum_{n^\alpha r^2\leq q\leq n^\beta r^2} \P(H_{(n+2)r}=q)\,
\P(H_{nr}\leq n^\beta r^2,\f_{nr,(n+2)r}\notin \Theta_{n,r}\,|\, H_{(n+2)r}=q).
\end{align}

Fix $q$ such that $n^\alpha r^2\leq q\leq n^\beta r^2$. We have
\begin{align}
\label{decompo2}
&\P(H_{nr}\leq n^\beta r^2,\,\f_{nr,(n+2)r}\notin \Theta_{n,r}\mid H_{(n+2)r}=q)\nonumber\\
&\qquad= \sum_{p\leq n^\beta r^2} \P(H_{nr}=p,\,\f_{nr,(n+2)r}\notin \Theta_{n,r}\mid H_{(n+2)r}=q)\nonumber\\
&\qquad= \sum_{p\leq n^\beta r^2} 
\frac{\varphi_{nr}(p)}{\varphi_{(n+2)r}(q)}\,\sum_{\f\in \F_{p,q,2r}\backslash \Theta_{n,r}} \prod_{v\in \f^*} \theta(c_v)
\end{align}
by formula \eqref{law-forest}.

For $p\leq n^\beta r^2$, we first bound the quantity 
$$\frac{\varphi_{nr}(p)}{\varphi_{(n+2)r}(q)}\leq M \,\frac{p}{q}\,\frac{(\pi_{nr})^{p-1}}{(\pi_{(n+2)r})^{q-1}}$$
where $M$ is a constant and the quantities $\pi_r$ were defined in Proposition \ref{hull-peri2}. Since 
$\pi_{nr}\leq 1$, we obtain that
\begin{equation}
\label{first-bd}
\frac{\varphi_{nr}(p)}{\varphi_{(n+2)r}(q)}\leq M \,\frac{p}{q}\,(\pi_{(n+2)r})^{-q}\leq M\,n^{\beta-\alpha}\,\exp( \frac{Aq}{n^2r^2})
\leq M\,n^{\beta-\alpha}\,\exp(A\,n^{\beta-2})
\end{equation}
with some constant $A$. On the other hand, the quantity
\begin{equation}
\label{bad-trees}
\sum_{p\leq n^\beta r^2} 
\sum_{\f\in \F_{p,q,2r}\backslash \Theta_{n,r}} \prod_{v\in \f^*} \theta(c_v)
\end{equation}
is bounded above by the probability that a forest 
of $q$ independent Galton-Watson trees with offspring distribution $\theta$
(truncated at level $2r$) is not in $\Theta_{n,r}$. For each tree in this forest, 
the probability that it satisfies property $(P)_r$
is at least $a_0/r^2$. The quantity \eqref{bad-trees} is thus bounded above by
$$\P(\ve_1+\cdots+\ve_q < n^\delta)$$
where the random variables $\ve_1,\ve_2,\ldots$ are i.i.d., with $\P(\ve_1=1)=1-\P(\ve_1=0)= a_0r^{-2}$. 
Since $q\geq n^\alpha r^2$ and $\alpha>\delta$, standard estimates on the
binomial distribution show that the quantity in the last display is bounded above
by $\exp(-\tilde a_0n^\alpha)$
for all $n$ sufficiently large and for every $r\geq 1$, with some other constant $\tilde a_0>0$.
Recalling \eqref{first-bd}, we get that the quantity \eqref{decompo2} is bounded above for $n$ large by
$$M\,n^{\beta-\alpha}\,\exp(A\,n^{\beta-2})\,\exp(-\tilde a_0n^\alpha).$$
Since $\beta-2<\alpha$, this shows that the left-hand side of \eqref{decompo} goes to 
$0$ faster than any negative power of $n$, uniformly in $r\geq 1$. 

By combining this observation with \eqref{first-esti} and \eqref{2-esti}, we obtain that,
for $n$ large enough,
\begin{equation}
\label{3-esti}
\P(\f_{nr,(n+2)r}\notin \Theta_{n,r})\leq C' n^{\frac{3}{2}(\alpha-2)}
\end{equation}
with a constant $C'$ independent of $r$. 

Let us argue on the event $\{\f_{nr,(n+2)r}\in \Theta_{n,r}\}$. If
$\tau$ is a tree of $\f_{nr,(n+2)r}$ with height $2r$, the vertices 
of $\tau$ at height $2r$ correspond to consecutive edges of $\cc_{nr}$, and if $v=v(\tau)$ is the 
last vertex in clockwise order that is incident to these edges,
we know
from Section \ref{sec:geodesics} that there is a downward geodesic path from a vertex $\tilde v$ of $\cc_{(n+2)r}$ (incident to the edge which is the root of $\tau$) to $v$, which has length exactly $2r$. Also, by the comments
of the end of Section \ref{dec-layers}, we know that $v(\tau)$ belongs to the boundary
of the standard hull of radius $nr$. Moreover, let $\tau$ and $\tau'$ be two distinct 
 trees of $\f_{nr,(n+2)r}$ with height $2r$, and assume that they both satisfy property $(P)_r$.
 Then the part of the boundary of the standard hull of radius $nr$ between $v(\tau)$
 and $v(\tau')$, in clockwise or in counterclockwise order, must contain at least $c_0r^2$ vertices: 
 This follows from the definition of property $(P)_r$ and the fact that, with each vertex $a$
 of $\tau$ (or of $\tau'$) at generation $2r-1$ having at least
 one child at generation $2r$ we can associate a vertex of the UIPQ --- namely, the 
 last vertex (in clockwise order) incident to the edges that are children of $a$ --- which belongs to
 the boundary of the standard hull of radius $nr$, as explained at the end of Section \ref{sec:trunc}. 

Write $\cc^\bullet_{nr}$ for the boundary of the standard hull $B^\bullet_{nr}(\pp)$
and $H^\bullet_{nr}$ for its perimeter (note that $H^\bullet_{nr}\leq 2 H_{nr}$). 
Also let $\pp^\bullet_{nr}$ stand for the complement of $B^\bullet_{nr}(\pp)$ 
in the UIPQ, viewed as an infinite quadrangulation with (simple) boundary $\cc^\bullet_{nr}$. 
By
the preceding observations,
the property
$\f_{nr,(n+2)r}\in \Theta_{n,r}$ implies that there are vertices 
$u_1,u_2,\ldots,u_k$ of $\cc^\bullet_{nr}$, with $k\geq n^\delta$, such that, for every $1\leq i\leq k$, 
$u_i$ is at graph distance $nr$ from the root vertex of the UIPQ and is connected to 
a vertex $\tilde u_i$ of $\cc_{(n+2)r}$ by a path of length $2r$, and moreover, if 
$i\not =j$, $u_i$ and $u_j$ are separated by at least $c_0r^2$ edges of $\cc^\bullet_{nr}$. 
Furthermore, write $\e_{n,r}$ for the event considered in Proposition \ref{main}:
$\e_{n,r}$ is the event where there exists a cycle
$\gamma$ of length smaller than $r$ that stays in $\cc(nr,(n+2)r)$, does not intersect $\cc_{(n+2)r}$,
and disconnects the root vertex of $\pp$ from infinity. On the event $\e_{n,r}\cap \{\f_{nr,(n+2)r}\in \Theta_{n,r}\}$,
for every $i\in\{1,\ldots,k\}$, the path from $u_i$ to $\tilde u_i$ must intersect the latter cycle: Indeed, if
we concatenate this path with a geodesic from the root vertex to $u_i$ and then with a path which goes
from the vertex $\tilde u_i$ to infinity and does not visit vertices at distance smaller than $(n+2)r$
from the root
(see Section \ref{sec:trunc-UIPQ}), we get a
path $\Gamma_i$ from the root vertex to infinity, which must intersect $\gamma$.
It follows that, for every $i,j\in\{1,\ldots,k\}$, we can construct a path of length at most $5r$ 
between $u_i$ and $u_j$ that stays
in $\cc(nr,(n+2)r)$. Since $u_i$ and $u_j$ both belong to $\cc^\bullet_{nr}$, a simple combinatorial argument shows that
this path can be required to stay in $\pp^\bullet_{nr}$. 

We then use the fact that, conditionally on $H^\bullet_{nr}$, 
$\pp^\bullet_{nr}$ is 
an 
infinite planar quadrangulation with a simple boundary of size $H^\bullet_{nr}$, which is independent of $B^\bullet_{nr}(\pp)$: This 
follows from the spatial Markov property of the UIPQ (we refer to
Theorem 5.1 in \cite{AS} for the UIPT, and the argument for the UIPQ is exactly the same). 
For every even integer $m\geq 2$, write $\pp^{(m)}$ for the UIPQ with simple boundary
of length $m$ (see \cite{CM}) and $\partial \pp^{(m)}$ for the collection of its boundary vertices.
Also denote the graph distance on the vertex set of $\pp^{(m)}$ by $\dg^{(m)}$. Let 
$\e^{(m,n,r)}_1$ stand for the event where there are at least $k=\lceil n^\delta\rceil$ vertices 
$v_1,\ldots,v_k$ of $\partial \pp^{(m)}$ such that, if $i\not =j$, $v_i$ and $v_j$
are separated by at least $c_0r^2$ edges of $\partial \pp^{(m)}$, and moreover
$$\dg^{(m)}(v_{i},v_{j})\leq 5r.$$
We will verify that $\sup\{\P(\e^{(m,n,r)}_1):m\leq 2n^\beta r^2\}$ decays exponentially in $n$
uniformly in $r$. Since by previous observations, we know that
$$\P(\e_{n,r}\cap \{\f_{nr,(n+2)r}\in \Theta_{n,r}\})\leq \sup\{\P(\e^{(m,n,r)}_1):m\leq 2n^\beta r^2\}$$
Proposition \ref{main} will follow from the bound \eqref{3-esti} (observe that we can choose $\alpha>0$
small so that $\frac{3}{2}(2-\alpha)$ is as close to $3$ as desired). 

In order to get the preceding exponential decay, we first replace the 
UIPQ with simple boundary $\pp^{(m)}$ by the UIPQ with general boundary
of the same size, which we denote by $\wt\pp^{(m)}$, and without risk of confusion, we
keep the notation $\dg^{(m)}$ for the graph distance (see again \cite{CM} for the definition of the 
UIPQ with general boundary). We
write 
$\mathbf{c}_0,\mathbf{c}_1,\ldots,\mathbf{c}_{m-1}$ for the ``exterior'' corners of the boundary of $\wt\pp^{(m)}$
enumerated in clockwise order starting from the root corner. Consider the event $\e^{(m,n,r)}_2$ where
one can find integers $0\leq p_1<p_2<\cdots<p_k<m$, with $k=\lceil n^\delta\rceil$, such that
$p_{i+1}-p_i\geq c_0r^2$ for every $1\leq i<k$, and $p_1+m-p_k \geq c_0r^2$,
and furthermore 
$$\dg^{(m)}(\mathbf{c}_{p_i},\mathbf{c}_{p_j})\leq 5r$$
whenever $i\not =j$. In order to bound the probability of $\e^{(m,n,r)}_2$, recall that,
from the results of \cite{CM} about infinite planar quadrangulations with
a boundary, we can assign labels $\ell(0),\ell(1),\ldots,\ell(m-1)$ to the corners $\mathbf{c}_0,\mathbf{c}_1,\ldots,\mathbf{c}_{m-1}$, which correspond to  ``renormalized'' distances from infinity, and are 
such that $\ell(0)=0$. Moreover, the sequence $(\ell(0),\ell(1),\ldots,\ell(m))$
(with $\ell(m)=0$) is a discrete bridge of length $m$, and we have,
for every 
$i,j\in\{0,1,\ldots,m-1\}$,
\begin{equation}
\label{cactus}
\dg^{(m)}(\mathbf{c}_i,\mathbf{c}_j)\geq 
\ell(i)+\ell(j) - 2 \,\max\Big( \min_{m\in[i,j]} \ell(m),\min_{m\in[j,i]} \ell(m)\Big)
\end{equation}
with the same convention for intervals as in Section \ref{bridge}.
The bound \eqref{cactus} follows from the ``treed bridge'' representation of $\wt\pp^{(m)}$
in \cite[Theorem 2]{CM} as  an instance of the so-called cactus bound (see \cite[Proposition 5.9(ii)]{LGM}
for the case of finite planar quadrangulations, and \cite[Lemma 3.7]{GM} for the case of 
finite quadrangulations with a boundary --- the proof in our infinite setting is exactly the same).
We can thus apply the bound of Lemma \ref{tec-lem} to the bridge $(\ell(0),\ell(1),\ldots,\ell(m))$
to get that $\P(\e^{(m,n,r)}_2)$ decays exponentially as $n\to\infty$ uniformly in
$m\leq 10n^{\beta}r^2$ and in $r\geq 1$.

We still need to verify that a similar exponential decay holds if we replace
$\P(\e^{(m,n,r)}_2)$ by $\P(\e^{(m,n,r)}_1)$. To this end, we rely on Theorem 4 of \cite{CM}, which states
that conditionally on the event where the size of its boundary is equal to $m'$, the
core of $\tilde\pp^{(m)}$ is distributed as $\pp^{(m')}$ (the core of $\tilde\pp^{(m)}$ is obtained
informally by removing the finite ``components'' of $\tilde\pp^{(m)}$ that can be disconnected 
from the infinite part by removing just one vertex, see Fig.2 in \cite{CM}). Furthermore, 
the probability that the boundary size of the core of 
$\tilde\pp^{(3m)}$ is equal to $m$ is bounded below by $c_1m^{-2/3}$ for some 
constant $c_1>0$, as shown
in the proof of Theorem 1 of \cite{CM}. Noting that the graph distance between
two vertices of the core only depends on the core itself (and not on the finite
components hanging off the core), we easily conclude
that
$$\P(\e^{(m,n,r)}_1)\leq (c_1m^{-2/3})^{-1}\,\P(\e^{(3m,n,r)}_2),$$
and this completes the proof of Proposition \ref{main}. \hfill$\square$

\medskip

Let us now explain how part (i) of Theorem \ref{kri-conj} is derived from Proposition \ref{main}.

\medskip
\noindent{\it Proof of Theorem \ref{kri-conj} (i).} Let $R\geq 1$ and let $\ve\in(0,1/2)$. Without loss of generality we can assume that
$\ve R\geq 1$. Suppose that $\gamma$ is a cycle  separating $B_R(\pp)$ from infinity, with length smaller
than $r:=\lfloor \ve R\rfloor$. Also let $n=\lfloor 1/\ve\rfloor$ so that $nr\leq R$. Then
the cycle $\gamma$ disconnects $B_{nr}(\pp)$ from infinity, which also
implies that it disconnects $B^\bullet_{nr}(\pp)$ from infinity
(in particular, $\gamma$ does not intersect $B^\bullet_{nr}(\pp)$).
Let $k\geq 0$ be the first integer such that $\gamma$ intersects the
annulus $\cc((n+k)r,(n+k+1)r)$. Then $\gamma$ is
contained in $\cc((n+k)r,(n+k+2)r)$ and does not intersect $\cc_{(n+k+2)r}$ (otherwise this would contradict the fact that $\gamma$
has length smaller than $r$). These considerations show that the event 
$\{L(R)\leq \ve R\}$ is contained in the union over $k\geq 0$ of the events where 
there exists a cycle of length smaller than $r$ that is contained in $\cc((n+k)r,(n+k+2)r)$,
does not intersect $\cc_{(n+k+2)r}$, and disconnects the root vertex
from infinity. Hence, if $\beta\in(1,3)$ is given, we deduce 
from Proposition \ref{main} that
$$\P(L(R)\leq \ve R) \leq \sum_{k=0}^\infty C'_\beta\,(n+k)^{-\beta} \leq \wt C_\beta\,n^{-\beta+1}=\wt C_\beta (\lfloor 1/\ve\rfloor)^{-\beta+1}.$$
The result of Theorem \ref{kri-conj} (i) follows.
\hfill$\square$

\section{Upper bound on the size of the separating cycle}
\label{sec:upper}

In this section, we prove Theorem \ref{kri-conj} (ii). 
We recall the notation $\f_{u,w}$ introduced before Corollary \ref{law-hull-annulus}, for integers $0\leq u<w$. We write $N_{u,w}$
for the number of trees of $\f_{u,w}$ that have maximal height $w-u$. As explained in
Section \ref{sec:geodesics}, for every integer $R\geq 3$,
one can find a cycle of the UIPT that disconnects the ball $B_{R-2}(\pp)$ from infinity and
whose length is bounded above by $2R\,N_{R,2R}$. In order to get bounds on $N_{R,2R}$, we
determine more generally the distribution of $N_{u,w}$ for any $1\leq u<w$.

\begin{proposition}
\label{tree-max}
The generating function of $N_{u,w}$ is given by
$$\E[a^{N_{u,w}}]= a\,\Bigg(\frac{9-\pi_w}{9-a\pi_w-(1-a)\pi_{w-u}}\Bigg)^{1/2}\,
\Bigg(\frac{1-\pi_w}{1-a\pi_w-(1-a)\pi_{w-u}}\Bigg)^{3/2}\;,$$
where we recall the notation
$$\pi_k=g^{(k)}_\theta(0)= 1- \frac{8}{(3+2k)^2 -1}.$$
Consequently, $N_{u,w}$ has the same distribution as $1 + U + V$, where 
$U$ and $V$ are independent,
$U$ follows the negative binomial distribution with parameters $(\frac{1}{2}, \frac{\pi_w-\pi_{w-u}}{9-\pi_{w-u}})$ and 
$V$ follows the negative binomial distribution with parameters $(\frac{3}{2}, \frac{\pi_w-\pi_{w-u}}{1-\pi_{w-u}})$. 
\end{proposition}

\proof Let $q\geq 1$ and first condition on the event $H_w=q$. By formula \eqref{law-forest},
$$\E[a^{N_{u,w}}\mid H_w=q]=
\sum_{\f\in \cup_{p\geq 1}\F_{w-u,p,q}} a^{N_{w-u}(\f)}\,\frac{\varphi_u(Y_{w-u}(\f))}{\varphi_w(q)}
\prod_{v\in\f^*} \theta(c_v),$$
where $Y_{w-u}(\f)$ is the number of vertices of the forest $\f$ at generation $w-u$, and 
$N_{w-u}(\f)$ denotes the number of trees of this forest with maximal height $w-u$. Using the
explicit formula \eqref{formula-phi} for $\varphi_u(p)$, we get
\begin{equation}
\label{tree-max-tec0}
\E[a^{N_{u,w}}\mid H_w=q]=\frac{f(u)}{\varphi_w(q)}\;\E_q\Big[Y_{w-u}\,\pi_u^{Y_{w-u}-1}\,a^{N_{w-u}}\Big],
\end{equation}
where 
$$f(u)=\frac{64}{3}\,\frac{3+2u}{((3+2u)^2-1)^2},$$
and under the probability measure $\P_q$, we consider a forest of $q$
independent Galton-Watson trees with offspring distribution $\theta$, 
$(Y_k)_{k\geq 1}$ denoting the associated Galton-Watson process. 

Let us compute the quantity $\E_q\Big[Y_{w-u}\,\pi_u^{Y_{w-u}-1}\,a^{N_{w-u}}\Big]$. 
To simplify notation, we set $m=w-u$. It is convenient
to write $\xi$ for a random variable distributed as the number of vertices at generation $m$
in a Galton-Watson tree with offspring distribution $\theta$ conditioned to be non-extinct 
at generation $m$. Then a simple argument shows that
\begin{equation}
\label{tree-max-tec1}
\E_q\Big[Y_{m}\,\pi_u^{Y_{m}-1}\,a^{N_{m}}\Big]
=\sum_{k=1}^\infty \P_q(N_{m}=k)\,k\,\E[\xi\pi_u^{\xi-1}]\,\Big(\E[\pi_u^\xi]\Big)^{k-1}\,a^k,
\end{equation}
and furthermore, 
$$
E[\pi_u^\xi]= \frac{g^{(m)}_\theta(\pi_u)-g^{(m)}_\theta(0)}{1-g^{(m)}_\theta(0)}=\frac{\pi_w-\pi_m}{1-\pi_m}\ ,\quad
E[\xi \pi_u^{\xi-1}]= \frac{\dot g^{(m)}_\theta(\pi_u)}{1-g^{(m)}_\theta(0)},$$
where $\dot g^{(m)}_\theta$ stand for the derivative of $g^{(m)}_\theta$.
We are thus led to the calculation of
$$\sum_{k=1}^\infty \P_q(N_{m}=k)\,k\,\Big(\frac{\pi_w-\pi_m}{1-\pi_m}\Big)^{k-1}\,a^k
= a\;\E_q\Big[N_{m}\,\Big(a\Big(\frac{\pi_w-\pi_m}{1-\pi_m}\Big)\Big)^{N_{m}-1}\Big].$$
Since $\E_q[a^{N_{m}}]=(a(1-\pi_m)+\pi_m)^q$ and so
$\E_q[N_{m}\,a^{N_{m}-1}]=q(1-\pi_m)(a(1-\pi_m)+\pi_m)^{q-1}$, we easily obtain that the
quantities in \eqref{tree-max-tec1} are equal to
$$q\,a\,\dot g^{(m)}_\theta(\pi_u)\;\Big(a(\pi_w-\pi_m)+\pi_m\Big)^{q-1}.$$

We substitute this in \eqref{tree-max-tec0} and then use the formula for the
distribution of $H_w$ in Proposition \ref{hull-peri2}. It follows that
$$\E[a^{N_{u,w}}]=
f(u)\,a\,\dot g^{(m)}_\theta(\pi_u)\; \sum_{q=1}^\infty 
\frac{K_w}{\varphi_w(q)}\,\kappa_q\,q\,\Big(a(\pi_w-\pi_m)+\pi_m\Big)^{q-1}\, \Big(2\pi_w\Big)^q.$$
To compute the sum of the series in the last display, first note that
$$q\,\frac{K_w}{\varphi_w(q)}\,\Big(2\pi_w\Big)^q=\frac{1}{\kappa_1}2^{q-1}.$$
It follows that
\begin{align*}
\E[a^{N_{u,w}}]&=\frac{f(u)\,a\,\dot g^{(m)}_\theta(\pi_u)}{\kappa_1}\;\sum_{q=1}^\infty \kappa_q
\,\Big(2(a(\pi_w-\pi_m)+\pi_m)\Big)^{q-1}\\
&=\frac{32\sqrt{3}}{\sqrt{\pi}}\,\frac{f(u)\,a\,\dot g^{(m)}_\theta(\pi_u)}{\kappa_1
\sqrt{(9-(a(\pi_w-\pi_m)+\pi_m))(1-(a(\pi_w-\pi_m)+\pi_m))^3}} 
\end{align*}
by \eqref{gen-C}. A simple calculation shows that
$$\dot g^{(m)}_\theta(\pi_u)= \frac{f(w)}{f(u)}.$$
Recalling that $\kappa_1=\frac{32}{\sqrt{3\pi}}$, we get
$$\E[a^{N_{u,w}}]=3 \,f(w)\,a\, \frac{1}{\sqrt{(9-(a(\pi_w-\pi_m)+\pi_m))(1-(a(\pi_w-\pi_m)+\pi_m))^3}} $$
and we just have to note that
$$f(w)=\frac{1}{3} \,\sqrt{(9-\pi_w)(1-\pi_w)^3},$$
by an easy calculation. \endproof

\noindent{\bf Remark.} Certain details of the previous proof could be simplified by using the interpretation suggested
in the remark following Corollary \ref{hull-peri}. 

\smallskip

\noindent {\it Proof of Part (ii) of Theorem \ref{kri-conj}}.  
As we already noticed, we have $L(R-2)\leq 2R\,N_{R,2R}$ for every $R\geq 3$, and we also note that $R\mapsto L(R)$
is nondecreasing. Then the proof boils down to verifying that
there exist constants $C'$ and $\lambda'>0$ such that, for every $R\geq 1$,
$\E[\exp(\lambda' N_{R,2R})] \leq C'$. Noting that
$$\frac{\pi_{2R}-\pi_R}{1-\pi_R}$$
is bounded above by a constant $\eta<1$, this follows from Proposition \ref{tree-max} and the form of 
the negative binomial distribution. \hfill$\square$

\smallskip

\noindent{\bf Remark.} We can also use Proposition \ref{tree-max} to get a lower bound on the
probability that is bounded above in Theorem \ref{kri-conj}(i). From Proposition \ref{tree-max},
we have immediately
$$\P(N_{u,w}=1)= 
\Bigg(\frac{9-\pi_w}{9-\pi_{w-u}}\Bigg)^{1/2}\,
\Bigg(\frac{1-\pi_w}{1-\pi_{w-u}}\Bigg)^{3/2}\;,$$
and taking $u=R$ and $w=\lfloor (1+\ve) R\rfloor$, we get
$$\P(N_{R,\lfloor (1+\ve) R\rfloor}=1) \build{\la}_{R\to\infty}^{} \Big(\frac{\ve}{1+\ve}\Big)^3.$$
Since we know that $L(R-2)\leq 2\ve R\,N_{R,\lfloor (1+\ve) R\rfloor}$, we conclude that
$$\liminf_{R\to\infty} \P(L(R)\leq 2\ve R) \geq \Big(\frac{\ve}{1+\ve}\Big)^3,$$
which should be compared with Theorem \ref{kri-conj}(i). Note that we get a lower bound by
(a constant times) $\ve^3$, whereas the upper bound of Theorem \ref{kri-conj}(i) gives $
\ve^\delta$ for every $\delta<2$. It is very plausible that one can refine the 
preceding argument, by considering all annuli of the
form $\cc_{\lfloor(1+n\ve)R\rfloor, \lfloor(1+(n+1)\ve)R\rfloor}$, and get also
a lower bound of order $\ve^2$ (of course one should deal with the lack of independence
of the variables $N_{\lfloor(1+n\ve)R\rfloor, \lfloor(1+(n+1)\ve)R\rfloor}$).

\section{Isoperimetric inequalities}
\label{sec:iso}

In this section, we prove Theorem \ref{iso1} and Proposition \ref{iso}.
We start with the proof of Proposition \ref{iso}, which is easier.

\medskip
\noindent{\it Proof of Proposition \ref{iso}.}
We first observe that we can assume that
$n$ is larger than some fixed integer, by then adjusting the constant $c_\ve$ if necessary.
In agreement with the notation of the introduction, write $|B^\bullet_r(\pp)|$ for the number of faces 
of the UIPQ
contained in the standard hull $B^\bullet_r(\pp)$.

By Theorem 5.1 in \cite{CLG0} (or \cite[Section 6.2]{CLG1}), we know that $r^{-4}|B^\bullet_r(\pp)|$  converges in
distribution as $r\to\infty$ to a limit which is finite a.s. It follows that we can
fix an integer $M>0$ such that, for every $r\geq 1$,
\begin{equation}
\label{iso-tec1}
\P(|B^\bullet_r(\pp)| < M\,r^4) \geq 1-\frac{\ve}{2}.
\end{equation}

On the other hand, for every $r\geq 1$ and every integer $N\geq 1$, let $E_{r,N}$ stand for the event where the minimal length of
a cycle separating $B_r(\pp)$ from infinity is greater than $r/N$. By Theorem \ref{kri-conj}, we can fix $N\geq 1$
large enough so that $\P(E_{r,N})>1-\ve/2$ for every $r\geq 1$. 

We fix $c>0$ such that $M(N+1)^4 c^4<1$. Then, let $n\geq 1$ large enough so that $cn^{1/4}\geq 1$
and $M((N+1)\lceil cn^{1/4}\rceil)^4<n$.
We argue on the event
$$\{|B^\bullet_{(N+1)\lceil cn^{1/4}\rceil}(\pp)|< n\} \cap E_{N\lceil c n^{1/4}\rceil, N}$$
which has probability at least $1-\ve$ by the preceding observations (we use \eqref{iso-tec1} and our choice of $c$).

We claim that, on latter event, we have $|\partial A| \geq c\,n^{1/4}$ for every $A\in\mathcal{K}$ such that $|A|\geq n$.
Indeed, writing $\Delta$ for the root vertex of $\pp$, we distinguish two cases:
\begin{itemize}
\item[$\bullet$] If $\dg(\Delta,\partial A)> N\lceil cn^{1/4}\rceil+1$, we note that the ball $B_{N\lceil cn^{1/4}\rceil}(\pp)$
is separated from infinity by the cycle $\partial A$, which implies that $|\partial A| \geq c\,n^{1/4}$ by the very definition
of the event $E_{N\lceil c n^{1/4}\rceil, N}$. 
\item[$\bullet$] If $\dg(\Delta,\partial A)\leq N\lceil cn^{1/4}\rceil+1$, then we argue by contradiction 
assuming that $|\partial A| < c\,n^{1/4}$ and in particular the diameter of $\partial A$ is bounded above by $\lceil c\,n^{1/4}\rceil-1$.
To simplify notation, we set $r_n=(N+1)\lceil cn^{1/4}\rceil$. The property $\dg(\Delta,\partial A)\leq N\lceil cn^{1/4}\rceil+1$ ensures that 
any vertex of $\partial A$ is at distance at most $r_n$ from $\Delta$, and therefore
any edge of $\partial A$ is incident to a vertex at distance at most $r_n-1$ from $\Delta$. It follows that any face incident to
an edge of $\partial A$ is contained in the hull $B^\bullet_{r_n}(\pp)$. Consequently, the whole boundary $\partial A$ is contained in $B^\bullet_{r_n}(\pp)$,
and so is the set $A$. 
In particular, $|A| \leq |B^\bullet_{r_n}(\pp)|<n$,
which is a contradiction with our assumption $|A|\geq n$. 
\end{itemize}
This completes the proof of Proposition \ref{iso}. \hfill$\square$

\medskip
Before we proceed to the proof of Theorem \ref{iso1}, we state a proposition which
is a key ingredient of this proof.

\begin{proposition}
\label{tail-hull}
There exists a constant $M_0$ such that, for every $r\geq 1$,
$$\E[|B^\bullet_r(\pp)|] \leq M_0\,r^4.$$
\end{proposition}

We postpone the proof of Proposition \ref{tail-hull} to the end of the section.

\medskip
\noindent{\it Proof of Theorem \ref{iso1}.} Thanks to Proposition \ref{tail-hull}
and Markov's inequality, we have
for every $r\geq 1$, for every $a>0$,
\begin{equation}
\label{markov-hull}
\P(|B^\bullet_r(\pp)|\geq a) \leq \frac{M_0r^4}{a}.
\end{equation}
We then proceed in a way very similar to the proof of Proposition \ref{iso}. Let $\delta\in(0,1/4)$ and
$\delta'=\delta/2$. For every integer $p\geq 1$, set
$$c_p=p^{-\frac{3}{4}-\delta}$$
and 
$$N_p=\lceil p^{\frac{1}{2}+\delta'}\rceil.$$
Recalling the notation $E_{r,N}$ in the proof of Proposition \ref{iso}, we first observe that
the bound of Theorem \ref{kri-conj} (i) implies
$$\sum_{p=1}^\infty \P(E_{N_p\lceil c_p2^{p/4}\rceil,N_p}^c) <\infty.$$
Similarly, the bound \eqref{markov-hull} gives
$$\sum_{p=1}^\infty \P(|B^\bullet_{(N_p+1)\lceil c_p 2^{p/4}\rceil}(\pp)| \geq 2^p)
<\infty.$$
The Borel-Cantelli lemma then shows that a.s. there exists an integer $p_0$
such that the event 
$$\{|B^\bullet_{(N_p+1)\lceil c_p 2^{p/4}\rceil}(\pp)| < 2^p\}\cap E_{N_p\lceil c_p2^{p/4}\rceil,N_p} $$
holds for every $p\geq p_0$. However, we can now reproduce the same arguments as in the end of the proof
of Proposition \ref{iso} (replacing $n$ by $2^p$, $N$ by $N_p$ and $c$ by $c_p$) to get that, if the latter event holds for some $p\geq 1$, then
for every $A\in\mathcal{K}$ such that $2^p\leq |A|\leq 2^{p+1}$, we have
$$|\partial A| \geq \lceil c_p2^{p/4}\rceil\geq 2^{-1/4} (\log\,2)\,|A|^{\frac{1}{4}}\,(\log|A|)^{-\frac{3}{4}-\delta}.$$
This completes the proof. \hfill$\square$

\medskip
We still have to prove Proposition \ref{tail-hull}. We note that the continuous analog of the UIPQ, 
or of the UIPT, is the Brownian plane introduced in \cite{plane}. For every 
$r>0$ the definition of the hull of radius $r$ makes sense in the Brownian plane
and the explicit distribution of the volume of the hull was computed in \cite{CLG0},
showing in particular that the expected volume is finite and (by scaling invariance) equal to a constant times $r^4$. On the other hand, in the case of the UIPT, M\'enard \cite{Men2} was recently able to compute the
exact distribution of the volume of hulls. Such an exact expression is not yet available in the
case of the UIPQ, and so we use a different method based on the skeleton decomposition.

Before we proceed to the proof of Proposition \ref{tail-hull}, we state a lemma 
concerning truncated quadrangulations with a Boltzmann distribution. Let $p\geq 1$.
We say that a random truncated quadrangulation $\mathbf{M}$ with boundary size $p\geq 1$
is Boltzmann distributed if, for every integer $n\geq 1$, for every 
$\mathcal{M}\in\Q^\tr_{n,p}$, $\P(\mathbf{M}=\mathcal{M})=Z(p)^{-1}\,12^{-n}$.

\begin{lemma}
\label{size-trunc}
There exists a constant $L_0>0$ such that, for every $p\geq 1$, if $\mathbf{M}_p$
is a Boltzmann distributed truncated quadrangulation with boundary size $p$,
$$\E[ \# \mathsf{Inn}(\mathbf{M}_p)] \build{\sim}_{p\to\infty}^{} L_0\,p^2.$$
\end{lemma}

\noindent{\bf Remark.} By analogy with the case of triangulations \cite[Proposition 9]{CLG1}, one expects that $p^{-2} \# \mathsf{Inn}(\mathbf{M}_p)$
converges in distribution to (a scaled version of) the distribution with density $(2\pi)^{-1/2}x^{-5/2}\,\exp(-\frac{1}{2x})$. 

\proof By definition,
$$\E[ \# \mathsf{Inn}(\mathbf{M}_p)]
=Z(p)^{-1}\,\sum_{n=p}^\infty n\,12^{-n}\,\#\Q^{\tr}_{n,p}.$$
In particular, \eqref{asym-trunc} shows that $\E[ \# \mathsf{Inn}(\mathbf{M}_p)]<\infty$.
From the definition of $U$ in Section \ref{sec:enu}, we have
$$Z(p)=\sum_{n=p}^\infty 12^{-n}\,\#\Q^{\tr}_{n,p}=[y^p]U(\frac{1}{12},y),$$
where as usual $[y^p]U(\frac{1}{12},y)$ denotes the coefficient of $y^p$ in the series expansion of $U(\frac{1}{12},y)$.
Then the explicit formula \eqref{formula-U}, and standard singularity analysis \cite[Corollary VI.1]{FS}, show that
\begin{equation}
\label{size-trunc1}
Z(p) \build{\sim}_{p\to\infty}^{} c' \,p^{-5/2}\,2^{-p}
\end{equation}
for some constant $c'>0$, whose exact value is
unimportant for our purposes. Similarly,
$$\sum_{n=p}^\infty n\,12^{-n}\,\#\Q^{\tr}_{n,p}
=[y^p] \frac{\partial U}{\partial x}(\frac{1}{12},y).$$
where the partial derivative is a left derivative at $x=1/12$. 
Formula (3) in \cite{Kr} gives
$$\frac{\partial U}{\partial x}(\frac{1}{12},y)
=-\frac{y^2}{2} - \frac{y}{2}\, \frac{y^2-10y-32}{\sqrt{(18-y)(2-y)}}.$$
and again singularity analysis leads to
\begin{equation}
\label{size-trunc2}
\sum_{n=p}^\infty n\,12^{-n}\,\#\Q^{\tr}_{n,p} \build{\sim}_{p\to\infty}^{} c'' \,p^{-1/2}\,2^{-p}
\end{equation}
for some other constant $c''>0$. The lemma now follows from
\eqref{size-trunc1} and \eqref{size-trunc2}. \endproof

\medskip
\noindent{\it Proof of Proposition \ref{tail-hull}.} We first observe that in the statement of the
proposition we may replace the standard hull $B^\bullet_r(\pp)$ by the truncated hull
$\h^\tr_r$. Indeed, the standard hull $B^\bullet_r(\pp)$ is contained in the truncated hull
$\h^\tr_{r+1}$. So let $N_{(r)}$ be the number of inner faces in the truncated hull
$\h^\tr_r$. We aim at proving that $\E[N_{(r)}]\leq M_0r^4$
for some constant $M_0$. 

Recall our notation $\f^\circ_{(r)}$ for skeleton of $\h^\tr_r$. As we already noticed in the proof of Corollary \ref{hull-peri}, 
$\f^\circ_{(r)}$ is distributed according to $\mu^\circ_{r,1}$. The fact that the distribution of 
$\h^\tr_{r}$ is $\Delta_{r,1}$ (Corollary \ref{law-hull}) yields that,
conditionally on $\f^\circ_{(r)}=\f$, the truncated quadrangulations $M_v$, $v\in \f^*$
associated with the ``slots'' are independent, and $M_v$ is Boltzmann distributed with
boundary size $c_v+1$, with the notation $c_v$ for the 
number of offspring of $v$ in $\f$. We then observe that, still on the event $\{\f^\circ_{(r)}=\f\}$,
$$N_{(r)}\leq \sum_{v\in \f^*} \#\mathsf{Inn}(M_v)$$
(see formula \eqref{formula-inner} in the proof of Lemma \ref{distri-hull}). Using Lemma \ref{size-trunc}, it follows that
$$\E[N_{(r)}\!\mid\! \f^\circ_{(r)}=\f] \leq L_0\,\E\Big[\sum_{v\in \f^*} (1+c_v)^2\Big].$$

As previously, it is convenient to use the notation $\f_{(r)}$ for the forest obtained 
by forgetting the distinguished vertex of $\f^\circ_{(r)}$ and applying a uniform
circular permutation to the trees of $\f^\circ_{(r)}$. From the last display, we
have also 
\begin{equation}
\label{tail-tec1}
\E[N_{(r)}] \leq L_0\,\E\Big[\sum_{v\in \f_{(r)}^*} (1+c_v)^2\Big],
\end{equation}
where we abuse notation by still writing $c_v$ for the number of offspring of the vertex $v$ of $\f_{(r)}$. 

In order to bound the expectation in the last display, we first consider vertices $v$
that are roots of trees in $\f^*_{(r)}$ (or equivalently which correspond to edges of $\cc_r$). 
In the forthcoming calculations, we also assume that $r\geq 2$. Let 
$c_{(1)},c_{(2)},\ldots,c_{(H_r)}$ denote the offspring numbers of the 
roots of the successive trees in $\f^*_{(r)}$. Using formula \eqref{law-forest} 
applied with $u=r-1$ and $w=r$, we get, for every $k\geq 1$,
$$\E\Big[\sum_{i=1}^{H_r} (1+c_{(i)})^2\, \Big|\,H_r=k\Big]
= \E\Bigg[ \frac{\varphi_{r-1}\Big(\xi_1+\cdots+\xi_k\Big)}{\varphi_r(k)}\,\sum_{i=1}^k (1+\xi_i)^2\Bigg],$$
where $\xi_1,\xi_2,\ldots$ are i.i.d. with distribution $\theta$. Recall formula
\eqref{law-peri} for the distribution of $H_r$, and the definition \eqref{formula-phi}
of $\varphi_r(k)$. Using also the estimate \eqref{asymp-C} for asymptotics of the
constants $\kappa_p$, we get, with some constant $L_1$,
$$\frac{\P(H_r=k)}{\varphi_r(k)} \leq \frac{L_1}{\sqrt{k}}.$$
We have therefore
$$\E\Big[\sum_{i=1}^{H_r} (1+c_{(i)})^2\,\mathbf{1}_{\{H_r=k\}}\Big]
\leq \frac{L_1}{\sqrt{k}}\, \E\Big[ \varphi_{r-1}\Big(\sum_{i=1}^k\xi_i\Big)\,\sum_{i=1}^k (1+\xi_i)^2\Big].$$
At this point, we again use \eqref{formula-phi} to see that there exist positive constants
$L_2,L_3,a_1$ such that, for every $\ell\geq 1$,
$$\varphi_{r-1}(\ell) \leq L_2\,\frac{\ell}{r^3}\,(1- \frac{a_1}{r^2})^\ell
\leq L_2\,\frac{\ell}{r^3}\,\exp(-\frac{a_1\ell}{r^2})\leq \frac{L_3}{r}\,\exp(-\frac{a_1\ell}{2r^2}).$$
It follows that
\begin{align*}
\E\Big[&\sum_{i=1}^{H_r} (1+c_{(i)})^2\,\mathbf{1}_{\{H_r=k\}}\Big]
\leq \frac{L_1L_3}{r\sqrt{k}}\,\E\Bigg[\Big(\sum_{i=1}^k (1+\xi_i)^2\Big) \exp\Big(-\frac{a_1}{2r^2}\sum_{i=1}^k\xi_i\Big)\Bigg]\\
&= \frac{L_1L_3\sqrt{k}}{r}\,\E\Big[(1+\xi_1)^2\,\exp(-\frac{a_1}{2r^2}\xi_1)\Big]\, \Big(
\E\Big[\exp(-\frac{a_1}{2r^2}\xi_1)\Big]\Big)^{k-1}.
\end{align*}
Using the asymptotics \eqref{equiv-theta} for $\theta(k)$, it is elementary to verify that
$$\E\Big[(1+\xi_1)^2\,\exp(-\frac{a_1}{2r^2}\xi_1)\Big]\leq L_4r$$
for some constant $L_4$. Moreover, we can also find a constant $a_2>0$ such that
$$\E\Big[\exp(-\frac{a_1}{2r^2}\xi_1)\Big]\leq 1-\frac{a_2}{r^2}.$$
We then conclude that
$$\E\Big[\sum_{i=1}^{H_r} (1+c_{(i)})^2\,\mathbf{1}_{\{H_r=k\}}\Big] \leq 
L_1L_2L_4\,\sqrt{k}\,(1-\frac{a_2}{r^2})^{k-1}.$$
By summing this estimate over $k\geq 1$, we get
$$\E\Big[\sum_{i=1}^{H_r} (1+c_{(i)})^2\Big] \leq L_5\,r^3$$
with some other constant $L_5$. 

A similar estimate holds if instead of summing over the roots of trees in the forest 
$\f_{(r)}$ we sum over vertices at generation $r-j$, for every $1\leq j\leq r-1$,
as this amounts to replacing $\f_{(r)}$ by the forest $\f_{(j)}$ (the case
$j=1$ requires a slightly different argument since we assumed $r\geq 2$
in the above calculations). By summing over $j$, recalling
\eqref{tail-tec1}, we conclude that
$$\E[N_{(r)}]\leq L_0\,\E\Bigg[\sum_{v\in \f_{(r)}^*} (1+c_v)^2\Bigg]\leq L_6 \sum_{j=1}^r j^3$$
with some constant $L_6$. This completes the proof of Proposition \ref{tail-hull}.

\section*{Appendix. Proof of Lemma \ref{tec-lem}}

We first note that $\pp_{k,K}(r)$ is trivially empty if $k>2K/(cr^2)$. If $k\leq 2K/(cr^2)$ and if we restrict our attention to
$k\leq k_0$ for some constant $k_0$, the bound of the lemma holds for any choice of $\gamma\in(0,1)$
by choosing the constant $C$ large enough. So we may assume that $k\geq k_0$ where $k_0$ can be taken large (but fixed).

Recall that $(b(0),b(1),\ldots,b(2K))$ stands for a discrete bridge of length $2K$. We first observe that,
for every $\ell \in\{0,1,\ldots,2K-1\}$, we can ``re-root'' the bridge $b(\cdot)$ at 
$\ell$ by setting, for every $j\in\{0,1,\ldots,2K\}$,
$$b_\ell(j):=\left\{\begin{array}{ll}
b(\ell+j)-b(\ell)&\hbox{if }\ell+j\leq 2K,\\
b(\ell+j-2K)-b(\ell)\quad&\hbox{if }\ell+j> 2K.
\end{array}
\right.
$$
Then $b_\ell(\cdot)$ is again a discrete bridge of length $2K$. Moreover the property defining 
the event $\pp_{k,K}(r)$ holds for $b(\cdot)$ with the sequence of times $(m_1,\ldots,m_k)$ if and only if it holds for $b_\ell(\cdot)$
with the sequence $(m^{(\ell)}_1,\ldots,m_k^{(\ell)})$ which is obtained by ordering the
representatives in $\{0,1,\ldots,2K-1\}$ of $m_1-\ell,\ldots,m_k-\ell$ modulo $2K$.

We start with a trivial observation. Let $0\leq m_1<m_2<\cdots<m_k<2K$ be integers. If $i_0\in\{1,\ldots,k\}$ is an index such that 
the minimal value of $b(\cdot)$ is attained in the interval $[m_{i_0},m_{i_0+1}]$ (there is at least one such value $i_0$), then,
for every $i,j\in\{1,\ldots,k\}$, $i\not =j$, the maximum
$$\max\Big( \min_{\ell\in[m_i,m_j]} b(\ell),\min_{\ell\in[m_j,m_i]} b(\ell)\Big)$$
is attained for the one among the two intervals $[m_i,m_j]$ and $[m_j,m_i]$ that does not contain $[m_{i_0},m_{i_0+1}]$.

Suppose that $0\leq m_1<m_2<\cdots<m_k<2K$ are such that the property of the definition of $\pp_{k,K}(r)$
holds, and that $k\geq 16$ is an integer multiple of $4$ (we can make the latter assumption without loss
of generality). As already mentioned, the property defining $\pp_{k,K}(r)$ still holds if 
$b(\cdot)$ is replaced by the re-rooted bridge $b_\ell(\cdot)$, for every $\ell\in\{0,1,\ldots,2K-1\}$, with the sequence 
$0\leq m^{(\ell)}_1<m^{(\ell)}_2<\cdots<m_k^{(\ell)}$ defined as explained above. Moreover, a  simple argument shows that there are at least 
$\frac{k}{2}$ values of $i\in\{1,\ldots,k\}$ such that, if $\ell\in(m_i,m_{i+1}]$, the minimum of
$b_\ell(\cdot)$ is attained in an interval $[m^{(\ell)}_j,m^{(\ell)}_{j+1}]$ with $\frac{k}{4}\leq j<\frac{3k}{4}$. Suppose that
$\ell$ is chosen uniformly at random in $\{0,1,\ldots,2K-1\}$, conditionally given $b(\cdot)$: the conditional probability for
the minimum of
$b_\ell(\cdot)$ to be attained in an interval $[m^{(\ell)}_j,m^{(\ell)}_{j+1}]$ with $\frac{k}{4}\leq j<\frac{3k}{4}$ is thus at least
$$\frac{\frac{k}{2} \times cr^2}{2K}.$$
It follows that
$$\frac{ckr^2}{4K} \,\P(\pp_{k,K}(r)) \leq \P(\pp^*_{k,K}(r))$$
where $\pp^*_{k,K}(r)$ is defined as $\pp_{k,K}(r)$ but imposing the additional constraint that the minimum of $b(\cdot)$
is attained in an interval $[m_{i_0},m_{i_0+1}]$ with $\frac{k}{4}\leq i_0<\frac{3k}{4}$. 

If $\pp^*_{k,K}(r)$ holds with the sequence $(m_1,\ldots,m_k)$, at least one of the two properties $m_{k/4}<K$ or $m_{3k/4}>K$
holds. We write $\pp^{**}_{k,K}(r)$ for the event where $\pp^*_{k,K}(r)$ holds and $m_{k/4}<K$ 
and we will bound the probability of $\pp^{**}_{k,K}(r)$ (the other case where $m_{3k/4}>K$ can be treated 
by time-reversal and leads to the same bound). 

Let us argue on the event $\pp^{**}_{k,K}(r)$. Using the definition of $\pp_{k,K}(r)$ and the trivial observation made at
the beginning of the proof, we note that, if $1\leq i\leq k/4$, 
$$b(m_k)+b(m_i)- 2\min_{\ell\in[m_k,m_i]} b(\ell) < 5r$$
and in particular
\begin{equation}
\label{tec-lem1}
b(m_i) < \min_{\ell\in[m_k,m_i]} b(\ell) + 5r\leq \un{b}(m_i) +5r,
\end{equation}
using the notation $\un{b}(\ell)=\min\{b(j):0\leq j\leq \ell\}$. 

We fix an integer $n\geq 1$ such that, if $S(0),S(1),\ldots$ is a simple random walk on $\Z$
started from $0$, the quantity
$$\P\Big(\min_{0\leq \ell\leq \lfloor ncr^2\rfloor} S(\ell) \geq -10r\Big)$$
is bounded above by a constant $\alpha<1$ independent of $r\geq 1$. Notice that the choice of
$n$ only depends on $c$. In what follows we assume that $k$ is large enough so that
$(\frac{k}{4}-1)/n \geq 1$ (recall the first observation of the proof).

We then define by induction $T_1=0$ and, for every integer $p\geq 1$,
$$T_{p+1}:=\inf\{\ell\geq T_p+ncr^2: b(\ell)\leq \un{b}(\ell)+5r\}\wedge K,
$$
where $\inf\varnothing=\infty$ as usual.
Recalling that we argue on $\pp^{**}_{k,K}(r)$, we notice that there are
at least $\lfloor(\frac{k}{4}-1)/n\rfloor$ (consecutive) values of $p\geq 1$ such that $m_1\leq T_p\leq m_{k/4}$. Indeed, the
first time $T_p$ that exceeds $m_1$ must be smaller than $m_{n+1}$ (by \eqref{tec-lem1} and our
assumption $m_{i+1}-m_i>cr^2$), the next one must be smaller than $m_{2n+1}$ if $2n+1\leq k/4$, and so on.
Moreover, if $p\geq 1$ is such that $m_1\leq T_p<T_{p+1}\leq m_{k/4}$, we have
$$
b(T_p)\leq \un{b}(T_p)+5r \leq b(m_1)+5r\leq \min_{\ell\in[m_1,m_{k/4}]} b(\ell) +10r\leq \min_{\ell\in[T_p,T_{p+1}]}b(\ell) +10r,$$
where in the third inequality we use the fact that $ b(m_1)\leq \min\{b(\ell) :\ell\in[m_1,m_{k/4}]\}+5r$, from the definition of
$\pp_{k,K}(r)$. Set $N_k=\lfloor(\frac{k}{4}-1)/n\rfloor-1$. We have obtained that 
$\pp^{**}_{k,K}(r)$ is contained in the event
\begin{align*}
e_{k,K}(r):=&\bigcup_{j=0}^\infty\Big\{ T_{j+N_k+1}<K\\
&\hbox{ and } b(T_{j+p})\leq \min_{\ell\in[T_{j+p},T_{j+p+1}]} b(\ell) + 10r,
\hbox{ for every }1\leq p\leq N_k\Big\}.
\end{align*}
Note that in the last display we can restrict the union to values of 
$j\in\{0,1,\ldots,\lfloor \frac{K}{ncr^2}\rfloor -1\}$, since by construction $T_p\geq (p-1)ncr^2$ if $T_p<K$. 

Recall that $S(0),S(1),\ldots$ is a simple random walk on $\Z$ started from $0$, and let 
$\wt T_1,\wt T_2,\ldots$ be the stopping times defined like $T_1,T_2,\ldots$
by replacing $(b(0),\ldots,b(2K))$ by $(S(0),\ldots,S(2K))$
and removing ``$\wedge \,K$''. 
We know that the distribution of $(b(0),b(1),\ldots,b(K))$ is absolutely continuous
with respect to that of $(S(0),S(1),\ldots,S(K))$, with a Radon-Nikodym
derivative that is bounded by a constant $M$ independent of $K$. It follows that
\begin{align*}
\P(\e_{k,K}(r))&\leq M\sum_{j=0}^{\lfloor\frac{K}{ncr^2}\rfloor -1} \!\P\Big(S(\wt T_{j+p})\leq \min_{\ell\in[\wt T_{j+p},\wt T_{j+p+1}]} S(\ell) + 10r,
\forall 1\leq p\leq N_k\Big)\\
&\leq M\times \lfloor\frac{K}{ncr^2}\rfloor\times \P\Big(\min_{0\leq \ell\leq ncr^2} S(\ell) \geq -10r\Big)^{N_k}\\
&\leq M\times \lfloor\frac{K}{ncr^2}\rfloor\times \alpha^{N_k}
\end{align*}
using the strong Markov property of $S$ in the second line, and our choice of $n$ in the last one. We conclude that
$$\P(\pp^{**}_{k,K}(r))\leq \P(\e_{k,K}(r))\leq M\times \lfloor\frac{K}{ncr^2}\rfloor\times \alpha^{N_k}$$
and since we have 
$$\P(\pp_{k,K}(r))\leq \Big(\frac{ckr^2}{4K}\Big)^{-1}\times 2\,\P(\pp^{**}_{k,K}(r))$$
we get the bound of the lemma. \endproof

\medskip
\noindent{\bf Acknowledgement.} The first author wishes to thank Nicolas Curien for several enlightening
discussions. We are also indebted to two referees for a very careful reading of the manuscript and 
several helpful suggestions.

\end{document}